\documentclass[reqno]{amsart}
\usepackage[mathscr]{eucal}
\usepackage{ulem,mathrsfs,amssymb,amsmath,calligra,fontenc,amsbsy,backref,
bbm,dsfont,graphicx}
\usepackage[utf8]{inputenc}
\usepackage[reversemp,
    paperwidth=210mm,
    paperheight=297mm,
    top={26mm},
    headheight={5.5pt},
    headsep={5.6mm},
    text={150mm,237mm},
    marginparsep=5mm,
    marginparwidth=12mm,
    bindingoffset=10mm,
    footskip=10mm]{geometry}
\DeclareMathAlphabet{\mathpzc}{OT1}{pzc}{m}{it}
\numberwithin{equation}{section}
\newcommand{\W}{{\normalfont\textsf{W}^{2}_\sharp}}
\newcommand{{\M}}{{\normalfont\textsf{L}^{2}_\sharp}}

\newcommand{\ts}{{\tilde{\sf s}}}

\newcommand{\C}{\mathbb C}
\newcommand{\bsfW}{\boldsymbol{\mathsf W}}

\newcommand{\bsfU}{\boldsymbol{\mathsf U}}
\newcommand{\bsff}{\boldsymbol{\mathsf f}}
\newcommand{\bsfF}{\boldsymbol{\mathsf F}}

\newcommand{\bsfu}{\boldsymbol{\mathsf u}}
\newcommand{\bsfz}{\boldsymbol{\mathsf z}}
\newcommand{\bsfw}{\boldsymbol{\mathsf w}}
\newcommand{\bsfe}{\boldsymbol{\mathsf e}}
\newcommand{\comp}{{\mathbb C}}
\newcommand{\comps}{{\mbox{\tiny $\comp$}}}

\newcommand{\bfsf}[1]{{\textbf{\textsf{#1}}}}

\newcommand{\sfw}{{\sf w}}

\newcommand{\sfp}{{\sf p}}
\newcommand{\sfv}{{\sf v}}

\newcommand{\Cdot}{\!\cdot\!}

\newcommand{\Div}{\mbox{\rm div}\,}

\newcommand{\curl}{\mbox{\rm curl}\,}

\newcommand{\IdS}{\Int{\partial\Omega}{}}

\newcommand{\Int}[2]{{\displaystyle \int_{ #1}^{ #2}}}
\newcommand{\Lim}[1]{{\displaystyle \lim_{ #1}}}

\newcommand{\Sup}[1]{{\displaystyle \sup_{#1}}}

\newcommand{\Sum}[2]{{\displaystyle \sum_{#1}^{#2}}}

\newcommand{\beea}{\begin{eqnarray}}
\newcommand{\eeea}{\end{eqnarray}}
\newcommand{\ms}{\medskip\smallskip}

\newcommand{\bfe}{{\mbox{\boldmath $e$}} }

\newcommand{\bfz}{{\mbox{\boldmath $z$}} }
\newcommand{\0}{{\mbox{\boldmath $0$}} }

\newcommand{\BF}{\begin{footnotesize}}
\newcommand{\EF}{\end{footnotesize}}
\newcommand{\ode}[2]{{\displaystyle \frac{\mbox{$d #1$}}{\mbox{$d #2$}}}}

\newcommand{\bi}{\begin{itemize}}
\newcommand{\ei}{\end{itemize}}
\newcommand{\ed}{\end{document}}
\newcommand{\be}{\begin{equation}}
\newcommand{\bestar}{\begin{equation*}}
\newcommand{\ba}{\begin{array}}
\newcommand{\ea}{\end{array}}
\newcommand{\eestar}{\end{equation*}}
\newcommand{\eeq}[1]{\label{eq:#1}\end{equation}}
\newcommand{\real}{\mathbb R}
\newcommand{\compl}{\mathbb C}

\newcommand{\nat}{{\mathbb N}}
\newcommand{\bfpsi}{\mbox{\boldmath $\psi$}}

\newcommand{\bfxi}{\mbox{\boldmath $\xi$}}
\newcommand{\bfphi}{\mbox{\boldmath $\varphi$}}

\newcommand{\bfv}{{\mbox{\boldmath $v$}} }
\newcommand{\bfu}{{\mbox{\boldmath $u$}} }

\newcommand{\bfw}{{\mbox{\boldmath $w$}} }
\newcommand{\bff}{{\mbox{\boldmath $f$}} }

\newcommand{\bfG}{{\mbox{\boldmath $G$}} }

\newcommand{\bfN}{{\mbox{\boldmath $N$}} }
\newcommand{\bfh}{{\mbox{\boldmath $h$}} }

\newcommand{\calc}{{\mathcal C}}
\newcommand{\cald}{{\mathcal D}}

\newcommand{\calf}{{\mathcal F}}

\newcommand{\calh}{{\mathcal H}}
\newcommand{\cali}{{\mathcal I}}

\newcommand{\calk}{{\mathcal K}}
\newcommand{\call}{{\mathcal L}}

\newcommand{\cals}{{\mathcal S}}

\newcommand{\calu}{{\mathcal U}}
\newcommand{\calv}{{\mathcal V}}

\newcommand{\caly}{{\mathcal Y}}
\newcommand{\calw}{{\mathcal W}}
\newcommand{\calz}{{\mathcal Z}}
\newcommand{\bfsigma}{\mbox{\boldmath $\sigma$}}

\newcommand{\bfeta}{\mbox{\boldmath $\eta$}}
\newcommand{\bfT}{{\mbox{\boldmath $T$}} }
\newcommand{\bfV}{{\mbox{\boldmath $V$}} }

\newcommand{\bfW}{{\mbox{\boldmath $W$}} }
\newcommand{\bfF}{{\mbox{\boldmath $F$}} }

\newcommand{\bfg}{{\mbox{\boldmath $g$}} }
\newcommand{\bfn}{{\mbox{\boldmath $n$}} }

\newcommand{\half}{\mbox{$\frac{1}{2}$}}

\def\Bbb R{\real}
\def\hat{\widehat}
\def\tilde{\widetilde}
\def\bar{\overline}

\newcommand{\bfchi}{\mbox{\boldmath $\chi$}}
\newcommand{\bfalpha}{\mbox{\boldmath $\alpha$}}

\newcommand{\bfgamma}{\mbox{\boldmath $\gamma$}}

\newcommand{\body}{\mathscr B}

\newcommand{\bfcalg}{\mbox{\boldmath ${\mathcal G}$}}

\newcommand{\bfcalf}{\mbox{\boldmath ${\mathcal F}$}}

\newcommand{\ED}{\end{description}}

\def\tag{\renewcommand{\theequation}}

\newtheorem{theo}{Theorem}
\newtheorem{prop}[theo]{Proposition}
\newtheorem{lem}[theo]{Lemma}
\newtheorem{cor}[theo]{Corollary}
\newtheorem{defi}{Definition}
\newtheorem{rem}{Remark}
\newcommand{\Bd}{\begin{defi}\begin{rm}}
\newcommand{\Ed}{\end{rm}\end{defi}}

\newcommand{\Br}{\begin{rem}\begin{it}}
\newcommand{\Er}{\end{it}\end{rem}}

\newcommand{\Bp}{\begin{prop}\begin{sl}}
\newcommand{\EP}[1]{\end{sl}\label{prop:#1}\end{prop}}

\newcommand{\Bt}{\begin{theo}\begin{sl}}
\newcommand{\Et}{\end{sl}\end{theo}}
\newcommand{\Bl}{\begin{lem}\begin{sl}}
\newcommand{\El}{\end{sl}\end{lem}}

\renewcommand{\eqref}[1]{{\rm (\ref{eq:#1})}}
\newcommand{\Bc}{\begin{cor}\begin{sl}}
\newcommand{\Ec}{\end{sl}\end{cor}}
\newcommand{\ET}[1]{\end{sl}\label{theo:#1}\end{theo}}
\newcommand{\EDD}[1]{\end{rm}\label{defi:#1}\end{defi}}
\newcommand{\EL}[1]{\end{sl}\label{lem:#1}\end{lem}}
\newcommand{\theoref}[1]{{\rm Theorem \ref{theo:#1}}}

\newcommand{\ER}[1]{\end{it}\label{rem:#1}\end{rem}}
\newcommand{\EC}[1]{\end{sl}\label{cor:#1}\end{cor}}
\newcommand{\remref}[1]{{\rm Remark \ref{rem:#1}}}
\newcommand{\cororef}[1]{{\rm Corollary \ref{cor:#1}}}
\newcommand{\lemmref}[1]{{\rm Lemma \ref{lem:#1}}}

\renewcommand{\i}{{\rm i}}

\begin{document}

\title[Hopf Bifurcation]{Flow-induced Oscillations via Hopf Bifurcation\\ in a Fluid-Solid Interaction Problem}
\author{Denis Bonheure}
\address{
	Département de Mathématique\\
	Université Libre de Bruxelles\\
	Boulevard du Triomphe 155\\
	1050 Brussels - Belgium
}
\email{denis.bonheure@ulb.be}
\author{
Giovanni P. Galdi}
\address{Department of Mechanical Engineering and Materials Science\\
University of Pittsburgh\\
Benedum Engineering Hall 607 \\
Pittsburgh, PA 15261 - USA}
\email{galdi@pitt.edu}
\author{Filippo Gazzola}
\address{
	Dipartimento di Matematica\\
	Politecnico di Milano\\
	Piazza Leonardo da Vinci 32\\
	20133 Milan - Italy}
 \email{filippo.gazzola@polimi.it}

\date{April 5, 2024}

\maketitle

\begin{abstract}
We furnish necessary and sufficient conditions for the occurrence of a Hopf bifurcation in a particularly significant fluid-structure problem, where a Navier-Stokes liquid interacts with a rigid body that is  subject to an undamped elastic restoring force. The motion of the coupled system  is driven by a uniform flow at spatial infinity, with constant dimensionless velocity $\lambda>0$. In particular,
if  the relevant linearized operator meets suitable spectral properties,  there exists a threshold $\lambda_o>0$ above which a bifurcating  time-periodic branch stems out of the branch of steady-state  solutions. The most remarkable feature of our result is that no restriction is imposed on the frequency $\omega$ of the bifurcating solution, which may thus coincide with one of the natural structural frequencies $\omega_{\sf n}$ of the body. Therefore, resonance cannot occur as a result of this bifurcation. However, when $\omega\to\omega_{\sf n}$, the amplitude of oscillations may become very large when the fluid density is negligible compared to the mass of the body.
To our knowledge, our result is the first {\it rigorous} investigation of the existence of a Hopf bifurcation  in a fluid-structure interaction problem.

{\bf Keywords. }{Navier-Stokes equations for incompressible viscous fluids, Fluid-solid interactions, Stability, Hopf Bifurcation}

{\bf Primary AMS code.} 76D05, 74F10, 35B32, 35Q35, 76D03, 35B10
 \end{abstract}

\tableofcontents
\renewcommand{\theequation}{\arabic{section}.\arabic{equation}}

\section{Introduction}\label{chap:intro}

The flow of a viscous fluid around structures is a  fundamental problem that lies at the heart of the broad research area of Fluid-Structure-Interactions (FSI). A main feature of this problem regards the study of the oscillations (vibrations) induced by the fluid on the structure, an astounding real-life example of Hopf bifurcation. Such oscillations may lead either to useful and pleasant motions, like ringing  wind chimes or Aeolian harps, or else destructive consequences, as damage or even collapse of the structure. In regards to the
latter, of particular significance is the phenomenon of forced oscillation of suspension bridges,
induced by the vortex shedding of the fluid (air), which reflects into an oscillatory regime of
the wake. When the frequency of the wake approaches one of the natural structural frequencies of the bridge, the so-called
lock-in \cite{diana}, a parametric resonance may occur and lead to a structural failure. A well known and infamous
example of this phenomenon is the  collapse of the Tacoma Narrows Bridge \cite{ammann,scott}.
\par
In view of its fundamental importance in many  practical situations, the problem of flow-induced oscillations of structures has received a rather large contribution by the engineering community. The list of the  most relevant articles is too long to be included here, and we direct the reader to the books \cite{Bev,Dyr,PPDL}, the review article \cite{Will} and the references therein. The structure model commonly employed is a rigid body (in most cases, a circular cylinder) subject to a linear restoring elastic force, while the fluid is modeled by the Navier-Stokes equations \cite{Bev}. In this regard, particularly remarkable is the detailed study by Blackburn and Henderson \cite{black} of the wake structures and flow dynamics associated with
simulated two-dimensional flows past a circular cylinder, either stationary or in simple harmonic cross-flow oscillation, in a range of  frequencies near the natural shedding frequency of the fixed cylinder.
\par
However, flow-induced oscillation problems have not yet received a similar  attention from the mathematical community. Actually, it is our  belief that a
{\it rigorous} mathematical approach to the problem could lead, on the one hand, to a deeper understanding of the underlying physics of the
phenomenon, suggesting possible solutions to improve stability of the structure, and, on the other hand, propose new interesting
and challenging questions that might  be of interest to both mathematicians and engineers.
\par
For this reason, a few years ago we began a rigorous and systematic investigation of flow-induced oscillations. Precisely, in \cite{BGG,BBGGP} we introduced the simple, but significant, model problem in which a two-dimensional rectangular structure is subjected to a unidirectional elastic restoring force, while immersed in
two-dimensional channel flow of a Navier-Stokes liquid, driven by a time-independent Poiseuille flow.
The main focus of these articles is rather fundamental and concerns the existence of possible equilibrium configurations of the structure, at least for "small" data. Further relevant properties of this model have been recently studied, such as explicit thresholds for the uniqueness of the equilibrium configuration \cite{GazP,GazzSp}, well-posedness of the related initial boundary value problem \cite{Patri}, large-time behavior \cite{BoHiPaSpe} and existence of a global attractor \cite{GPP}. It should be underlined that the model used in all the works just cited presents two simplifying characteristics. First of all, it is two-dimensional, thus limiting the number of phenomena that can occur. Furthermore, it includes boundary walls whose presence can influence the interaction of the fluid with the structure, which is in fact the only fundamental question we are interested in investigating flow-induced vibrations.
\par
Due to this, very recently \cite{BoGaGa1}, we introduced a different model --inspired by \cite{Bev}-- where the structure, $\mathscr S$, is a three-dimensional sufficiently smooth body, $\mathscr B$, of arbitrary shape subject to a (possibly anisotropic) undamped linear\footnote{See, however, \remref{0.0}.} restoring force, and the Navier-Stokes liquid, $\mathscr L$, fills the whole space outside $\mathscr B$.  The motion of the coupled system, $\mathscr S$+$\mathscr L$, is driven by a uniform flow at spatial infinity, characterized by a  constant dimensionless velocity $\lambda>0$. In \cite{BoGaGa1} we carried out a preliminary study devoted to the existence, stability and steady-state bifurcation of  equilibria, where $\mathscr B$ is in a fixed configuration and $\mathscr L$ performs a steady-state flow. In particular, we showed the existence of some ``critical" $\lambda_s>0$ such that if $\lambda<\lambda_s$ the equilibrium is unique, stable and no oscillatory motion can occur.
\par
The objective of this article is to continue the investigation begun in \cite{BoGaGa1} and get to the heart of the issue of flow-induced vibrations, i.e., the relationship between the onset of time-periodic (Hopf) bifurcation and the resonance phenomenon. To our knowledge, ours is the first study addressing this type of problems in a rigorous fashion.
\par
We shall now describe  our main results along with the corresponding difficulties in proving them. Let $\lambda_o>0$ and denote by $(\lambda,{\sf s}(\lambda))$, $\lambda\in U(\lambda_o)$ a family of smooth equilibria in the neighborhood $U$ of $\lambda_o$. We thus aim at finding necessary and sufficient conditions on $(\lambda_o,{\sf s}(\lambda_o))$ ensuring that a time-periodic branch does exist in a neighborhood of $(\lambda_o,{\sf s}(\lambda_o))$. This question is investigated in the framework of a general time-periodic bifurcation theory introduced in \cite{GaBif}, which, for completeness,  we recall in Section \ref{sec:bif}, with full details presented in the Appendix. Unlike classical approaches \cite{CR1}, \cite[Theorem 79.F]{Z1}, ours  is able to handle flows in unbounded regions and, in particular, overcome the notorious issue of 0 being in the essential spectrum of the linearized operator \cite{Bab0,Bab,FN,Saz}. The pivotal point of our method is to find functional frameworks (possibly distinct)  where the linearized operators around $(\lambda_o,{\sf s}(\lambda_o))$ for the averaged (over a period) and ``purely oscillatory" components of the relevant fields are both Fredholm of index 0. Searching for the ``right" functional setting is exactly what
makes the application of the main theorem of \cite{GaBif} to the case in question anything but simple. Such a study  is carried out in Sections \ref{sec:spectrum} and \ref{sec:time_per}. The {\it crucial} result established in Section \ref{sec:time_per} (see \theoref{3.1}) is that the ``purely oscillatory" linear operator is Fredholm of index 0,   {\it whatever} the (finite) value of the material constants.
In other words, the Fredholm property holds also in case of ``resonance", namely, when the frequency of the oscillation coincides with the natural frequency of $\mathscr B$ or a multiple of it (for suspension bridges, this is the so-called lock-in effect \cite{diana}, leading to vortex induced vibrations of the structure).
We also show that this important property no longer holds in the limit situation when the ratio of the density of $\mathscr L$ to that of $\mathscr B$ tends to 0, as somehow expected on physical ground.
That the model here employed does not seem to be well suited to predict resonance is also confirmed by the recent work \cite{BoGa0} where it is proved that, whatever its frequency, a time-periodic velocity field $\bfV$ far from the body always generates at least one periodic oscillation.\par
With these preparatory results in hand, we are able to reformulate the bifurcation problem in a suitable functional setting (Section \ref{sec:rif}), and, with the help of the theory given in \cite{GaBif}, to show our main finding, collected in Section \ref{sec:Bif} (see \theoref{8.1}). In a nutshell, the latter states that, under certain conditions ensuring the existence of an analytic steady-state branch in the neighborhood of $(\lambda_0,{\sf s}(\lambda_0))$ proved in \cite[Lemma 15]{BoGaGa1},   there is an analytic family of time-periodic solutions emanating from $(\lambda_0,{\sf s}(\lambda_0))$, provided the linearized operator at $(\lambda_0,{\sf s}(\lambda_0))$ has a ``non-resonant" simple, purely imaginary eigenvalue crossing the imaginary axis at non-zero speed, when $\lambda$ crosses $\lambda_0$. Also, as expected, the bifurcation is supercritical. It is worth emphasizing that, because of the range of validity of the Fredholm property mentioned above, no restriction is imposed on the frequency of the bifurcating branch, thus ruling out the occurrence of  ``disruptive" resonance.
\par
As a final comment, we recall that,  due to the intrinsic complexity of the nondegeneracy conditions arising in the mathematical treatment of bifurcation, our analytic findings must be complemented or supported by numerical tests. For instance, the numerical computation in \cite{FGQ,RR} illustrate the required nondegeneracy conditions given in \cite{GaArma}. We also mention a recent research line based on fully rigorous computer assisted proofs \cite{argako,arioli}.\par
Of course, the analysis performed here  is by no means exhaustive and leaves out a number of fundamental questions. First and foremost, the stability properties of the bifurcating solution. On physical ground, it is expected that the time-periodic bifurcating branch becomes stable, while the equilibrium loses its stability. However, with our current knowledge, a rigorous proof of this assertion remains out of our reach and will be the object of future thoughts.
\par
The article is divided into several sections, each one regarding different aspects of the problem. Specifically, in Section~\ref{PrFo}  we present the relevant equations and furnish  the mathematical formulation of the problem. Moreover, after recalling some notation, we introduce the appropriate function spaces and collect some of their important characteristics in Section \ref{spaces}.
In Section \ref{sec:bif} we set up an abstract approach to time-periodic bifurcation for general dynamical systems.
Section~\ref{chap:bifu} is devoted to finding necessary and sufficient conditions for time-periodic bifurcation.  As previously mentioned, this is done via  the general abstract approach introduced in \cite{GaBif} and presented here in full details.  This approach requires, among other things, a careful analysis of the properties of the relevant linear operators. Specifically, in Subsection~\ref{sec:spectrum} we show that, given $\lambda_o>0$, the intersection of the spectrum of the (time-independent) linearization at $(\lambda_o,{\sf s}(\lambda_o))$ with the imaginary axis is constituted (at most) by a bounded sequence of eigenvalues of finite algebraic multiplicity that can cluster only at 0. Such a result is crucial, in that it makes plausible the basic requirement for a Hopf-like bifurcation, namely, the existence of an imaginary eigenvalue of algebraic multiplicity 1. In Subsection~\ref{sec:time_per} we show that the linearized ``purely oscillatory" operator, suitably defined, is Fredholm of index 0, whatever the value of the physical parameters. With all the above results in hand, in Subsection~\ref{sec:rif} we then formulate the bifurcation problem in a functional setting that fits the requirement of the general theory  previously presented. In the final Subsection~\ref{sec:Bif}, we prove necessary and sufficient conditions for bifurcation, with the latter ensuring the existence of a (subcritical or supercritical) family of time-periodic solutions in the neighborhood of $(\lambda_o,{\sf s}(\lambda_o))$, stemming out of the analytic branch of equilibria.

\section{Formulation of the Problem}\label{PrFo}
Let $\mathscr B$ be a rigid body moving in a Navier-Stokes liquid that fills the region $\Omega\subset\mathbb{R}^3$ outside $\mathscr B$ and whose flow becomes uniform at ``large" distances from $\mathscr B$, characterized by a constant velocity $\bfV\in\mathbb{R}^3$. We denote by $\Omega_0$ the  volume occupied by $\mathscr B$, which we assume throughout  to be the closure of a bounded domain of class $C^2$. An elastic restoring force $\bfF$ acts on $\mathscr B$, applied to its center of mass $G$, while a suitable active couple prevents its rotation. Therefore, the motion of $\mathscr B$ is translatory. In this situation, the  governing equations of   motion of the coupled system body-liquid when  referred to a body-fixed frame $\calf\equiv\{G,\bfe_i\}$  are given by \cite[Section 1]{Gah}
\be\ba{cc}\medskip\left.\ba{r}\medskip
\partial_t\bfv-\nu\Delta\bfv+\nabla p+(\bfv-{\bfgamma})\cdot\nabla\bfv=0\\
\Div\bfv=0\ea\right\}\ \ \mbox{in $\Omega\times(0,\infty)$}\,,\\ \medskip
\bfv(x,t)={\bfgamma}(t)\,, \ \mbox{ $(x,t)\in\partial\Omega\times(0,\infty)$}\,;\ \
\Lim{|x|\to\infty}\bfv(x,t)=\bfV\,,\ t\in(0,\infty)\,,\\
M\dot{\bfgamma}+\rho\Int{\partial\Omega}{} \mathbb T(\bfv,p)\cdot\bfn=\bfF \ \ \mbox{in $(0,\infty)$}\,.
\ea
\eeq{01}
In \eqref{01}, $\bfv$ and $p$ represent velocity and pressure fields of the liquid,  $\rho$ and $\nu$ its density and kinematic viscosity,  while $M$ and $\bfgamma=\bfgamma(t)$ denote mass of $\mathscr B$  and velocity of $G$, respectively. Moreover, we consider the Cauchy stress tensor
$$
\mathbb T_\nu(\bfz,\psi):=2\nu\,\mathbb D(\bfz)-\psi\,\mathbb I\,,\ \ \ \mathbb D(\bfz):=\half\left(\nabla\bfz+(\nabla\bfz)^\top\right)\,,
$$
where $\mathbb I$ is the $3\times3$ identity matrix and $\bfn$ is the unit outer normal at $\partial\Omega$, i.e. directed inside $\body$.\par
We assume that $\bfF$ depends linearly on the displacement $\bfchi(t):=\int\bfgamma(s){\rm d}s=\vec{GO}$, with $O$ fixed point,
namely
\be
\bfF(t)=-\mathbb B\cdot\bfchi(t),\quad t \ge0,
\eeq{ElFo}
where $\mathbb B$ is a symmetric, positive definite matrix (stiffness matrix).
Without loss of generality we take $\bfV=-V\bfe_1$, $V>0$.

\Br
The choice of a linear restoring force in the  constitutive equation \eqref{ElFo} is made {\it just} for simplicity of presentation. A careful reading of the proofs shows that one can deal similarly with a nonlinear restoring force of the type $\bfF=\mathbb B\cdot\bfchi+\bfg(\bfchi)$, where $\bfg(\bfchi)$ is sufficiently smooth and $|\bfg(\bfchi)|=o(|\bfchi|)$ as $|\bfchi|\to 0$.
\ER{0.0}

Writing all the involved quantities in a non-dimensional form, and setting $\bfu:=\bfv+\bfe_1$,
we may rewrite \eqref{01} in the following form
\be\ba{cc}\medskip\left.\ba{r}\medskip
\partial_t\bfu-\Delta\bfu+\nabla p=\lambda\,[\partial_1\bfu+(\dot{\bfchi}-\bfu)\cdot\nabla\bfu]\\
\Div\bfu=0\ea\right\}\ \ \mbox{in $\Omega\times(0,\infty)$}\,,\\ \medskip
\bfu(x,t)={\dot{\bfchi}}(t)+\bfe_1\,, \ \mbox{ $(x,t)\in\partial\Omega\times(0,\infty)$}\,;\ \
\Lim{|x|\to\infty}\bfu(x,t)=\0\,,\ t\in(0,\infty)\,,\\
\ddot{\bfchi}+\mathbb{A}\cdot \bfchi+\varpi\Int{\partial\Omega}{} \mathbb T_1(\bfu,p)\cdot\bfn=\0 \ \ \mbox{in $(0,\infty)$}\,,
\ea
\eeq{02}
with
$$
\mathbb{A}:=\frac{L^4}{M\nu^2}\mathbb B\,,\    \ \varpi:=\frac{\rho L^3}{M}\,,\ \ \lambda:=\frac{VL}{\nu}\,,
$$
In the sequel, we just write  $\mathbb T$ instead of $\mathbb T_1$. Since $\mathbb{A}$ is positive definite, $\mathbb{A} $ is invertible and there exists $\sf a, b>0$ such that for every $\bfalpha\in \mathbb C^3$,
$$\sf a\, \|\bfalpha\|^2\le \bfalpha^*\cdot\mathbb{A} \cdot\bfalpha\le \sf b\, \|\bfalpha\|^2.$$
As it will be clear in the proofs, when we indicate that a constant in the estimates depends on $\mathbb{A}$, the dependence is obviously limited to the extremal eigenvalues of $\mathbb{A}$ (as it is symmetric), i.e. either on $\sf a$ or $\sf b$, or both.\par
We are interested in the existence of  solutions $(\bfu,p,\bfchi)$ to \eqref{02}  bifurcating from a branch of equilibrium configurations. In order to make this statement more precise,  let ${\sf s}_0=(\bfu_0,p_0,\bfchi_0)$ be an equilibrium solution to \eqref{02} corresponding to a given $\lambda$, namely,
\be\ba{cc}\medskip\left.\ba{r}\medskip
-\Delta\bfu_0+\nabla p_0=\lambda\,(\partial_1\bfu_0-\bfu_0\cdot\nabla\bfu_0)\\
\Div\bfu_0=0\ea\right\}\ \ \mbox{in $\Omega$}\,,\\ \medskip
\bfu_0(x)=\bfe_1\,, \ \mbox{ $x\in\partial\Omega$}\,;\ \
\Lim{|x|\to\infty}\bfu_0(x)=\0\,,\\
\mathbb{A} \cdot \bfchi_0+\varpi\Int{\partial\Omega}{} \mathbb T(\bfu_0,p_0)\cdot\bfn=\0\,.
\ea
\eeq{03}
From the physical viewpoint, $\bfchi_0$ represents the (non-dimensional and rescaled) elongation of the spring necessary to keep $\mathscr B$ in place.\par
We look for a value $\lambda_o> 0$  at which a time-periodic motion
may branch out of ${\sf s}_0(\lambda_o)$. Writing
$$
\bfu=\bfu_0(\lambda)+\bfw(t;\lambda)\,,\ p=p_0(\lambda)+{\sfp}(t;\lambda)\,,\ \ \bfchi=\bfchi_0(\lambda)+ \bfeta(t;\lambda)\,,
$$
equations \eqref{02} yields the first-order in time system
\be\ba{cc}\medskip\left.\ba{lr}\medskip
& \partial_t\bfw-\Delta\bfw+\nabla {\sfp}   \medskip  =\lambda\,[\partial_1\bfw-\bfu_0\Cdot\nabla\bfw+
({\bfsigma}-\bfw)\Cdot\nabla\bfu_0 -({\bfsigma}-\bfw)\Cdot\nabla\bfw]\\
& \Div\bfw=0\ea\right\}\ \ \mbox{in $\Omega\times(0,\infty)$}\,,\\ \medskip
\bfw(x,t)={{\bfsigma}}(t)\,, \ \mbox{ $(x,t)\in\partial\Omega\times(0,\infty)$}\,;\ \
\Lim{|x|\to\infty}\bfw(x,t)=\0\,,\ t\in(0,\infty)\,,\\
\dot{\bfsigma}+\mathbb{A} \cdot \bfeta+\varpi\Int{\partial\Omega}{} \mathbb T(\bfw,{\sfp})\Cdot\bfn=\0\,;\ \ \dot{\bfeta}(t)=\bfsigma(t) \ \ \mbox{in $(0,\infty)$}\,,
\ea
\eeq{04}
where $(\bfw,{\sfp},\bfeta,\bfsigma)\in\real^3\times\real\times\real^3\times\real^3$ are the four unknowns.

Our goal is to determine sufficient conditions on ${\sf s}_0(\lambda_o)$, for the existence
of a non-trivial family of time-periodic solutions to \eqref{04}, $(\bfw(\lambda), {\sfp}(\lambda),\bfeta(\lambda),\bfsigma(\lambda))$,
$\lambda\in U(\lambda_o)$, of (unknown) frequency $\zeta =\zeta (\lambda)$,  such that $(\bfw(\lambda), {\sfp}(\lambda),\bfeta(\lambda),\bfsigma(\lambda))\to (\0,\0,\0,\0)$ as $\lambda\to\lambda_o$. We will achieve this objective by employing the general approach to time-periodic bifurcation presented in Section~\ref{sec:bif}.
The implementation of this approach to the case at hand requires a detailed study of the properties of several linearized problems associated to \eqref{04}, which will be established in Section~\ref{chap:bifu}.

\section{Functional Background}\label{spaces}
This section is dedicated to introduce the relevant functional spaces we need.
With the origin of coordinates in the interior of $\Omega_0$, we set
$$
B_R:=\{x\in\real^3:\,|x|<R\},\ \Omega_R:=\Omega\cap B_R,\
\Omega^R:=\Omega\backslash\bar{\Omega_R}\quad\forall R>R_*:={\rm diam}\, \Omega_0.$$
For a domain $A\subset\real^3$, we use the following standard notations :
\begin{itemize}
\item $L^q=L^q(A)$ is the Lebesgue space with norm $\|\cdot\|_{q,A}$ and scalar product $(\cdot\ ,\,\cdot )_A$  for $q=2$;
\item $W^{m,2}=W^{m,2}(A)$, $m\in\nat$ is the classical Sobolev space of order $m$ with norm $\|\cdot\|_{m,2,A}$ ;
\item $D^{m,q}=D^{m,q}(A)$ is the homogeneous Sobolev space with semi-norm $\sum_{|l|=m}\|D^lu\|_{q,A}$ ;
\item $D_0^{1,2}=D_0^{1,2}(A)$ is the completion of $C_0^\infty(A)$ in the norm $\|\nabla(\cdot)\|_{2,A}$.
\end{itemize}
We will usually omit the subscript ``$A$" in the norms, unless confusion arises.

\medbreak

We first recall the relevant function spaces for the steady part of the solution, see also \cite{BoGaGa1}. We start with
\[
\calk=\mathcal K(\real^3):=\big\{\bfphi\in C_0^\infty({\real^3}):
\exists\, \hat{\bfphi}\in\real^3 \text{ s.t. }\bfphi(x)\equiv\hat{\bfphi} \text{ in a neighborhood of }\Omega_0\big\}
\]
that we endow with the scalar product
\be
\langle \bfphi,\bfpsi\rangle:=
\varpi^{-1}\,\hat{\bfphi}\cdot\hat{\bfpsi}+(\bfphi,\bfpsi)_\Omega\,,\quad\forall\bfphi,\bfpsi\in\calk\,.
\eeq{0.0}
This scalar product induces the norm $\|\cdot\|_{\calk}$. We then introduce
\begin{itemize}
\item $\calc=\mathcal C(\real^3):=\{\bfphi\in\calk(\real^3):\ \Div\bfphi=0\ \mbox{in $\real^3$}\}\,;$
\item $\calc_0=\mathcal C_0(\Omega):=\{\bfphi\in\calc(\real^3): \hat{\bfphi}=\0\}\,.$
\end{itemize}
The closure of $\calk(\real^3)$ and $\calc(\real^3)$ with respect to the norm $\|\cdot\|_{\calk}$ will be denoted by
$$\mathcal L^2=\call^2(\real^3):= \overline{\calk(\real^3)}^{\, \|\cdot\|_{\calk}}\, \text{ and } \quad \mathcal H=\calh(\real^3):= \overline{\calc(\real^3)}^{\, \|\cdot\|_{\calk}}\,,$$
whereas the closure of $\calc(\real^3)$ and $\calc_0(\Omega)$ with respect to the norm $\|\mathbb D(\cdot)\|_2$ will be denoted by
$$
\cald^{1,2}=\cald^{1,2}(\real^3):=\overline{\calc(\real^3)}^{\, \|\mathbb D(\cdot)\|_2}\, \text{ and } \quad \cald_0^{1,2}=\cald_0^{1,2}(\Omega):=\overline{\calc_0(\Omega)}^{\, \|\mathbb D(\cdot)\|_2}.$$
With these notations at hand, we set
\be
Z^{2,2}:=W^{2,2}(\Omega)\cap \cald^{1,2}(\real^3)\,.
\eeq{zetaa}
Obviously, $\cald_0^{1,2}(\Omega)\subset\cald^{1,2}(\real^3)$.
Denoting the dual space of $\cald_0^{1,2}(\Omega)$ by $\cald_0^{-1,2}(\Omega)$, endowed with the norm
$$
|\bff|_{-1,2}=\Sup{\mbox{\footnotesize $\ba{c}\bfphi\in \calc_0(\Omega)\\ \|\nabla\bfphi\|_2=1\ea $}}\left|\left(\bff,\bfphi\right)\right|\,,
$$
we define the spaces
\be
\caly:=\cald_0^{-1,2}(\Omega)\cap \calh(\real^3)\,,\ \ Y:=\cald_0^{-1,2}(\Omega)\cap \call^2(\real^3)\,,
\eeq{ips}
with associated norms
$$
\|\bfg\|_{\caly}=\|\bfg\|_{Y}:=\|\bfg\|_{2}+|\bfg|_{-1,2}+|\hat{\bfg}|\,,
$$
and in turn
\be
\begin{split}
X=X(\Omega):=\{\bfu\in\cald^{1,2}_0(\Omega):\, \partial_1\bfu\in \cald_0^{-1,2}(\Omega)\},\\ X^2=X^2(\Omega):=\left\{\bfu\in X(\Omega): D^2\bfu\in L^2(\Omega)\right\}\,.
\end{split}
\eeq{ics}
Both $X$ and $X^2$ are (reflexive, separable) Banach spaces when equipped with the norms\footnote{the norms $\|\nabla(\cdot)\|_2$ and $\|\mathbb D(\cdot)\|_2$ are equivalent in $\cald^{1,2}_0$, see e.g. \cite[Lemma 3]{BoGaGa1}.}
$$
\|\bfu\|_X:=\|\nabla\bfu\|_2+|\partial_1\bfu|_{-1,2}\,,\ \
\|\bfu\|_{X^2}:=\|\bfu\|_{X}+\|D^2\bfu\|_{2}\,,
$$
see for instance \cite[Proposition 65]{GaCe}.

We next introduce spaces of time-periodic functions.
A function $\bfw:\Omega\times \real\mapsto \real^3$ is
{\it $2\pi$-periodic}, if for a.e.\ $t\in \real$, $\bfw(\cdot,t+~2\pi)=\bfw(\cdot,t)$,
and we use the standard notation for its average
\be
{\bar \bfw(\cdot)}:=\frac{1}{2\pi}\int_{0}^{2\pi}\bfw(\cdot,t){\rm d}t\,,
\eeq{media}
whenever the integral is defined.
If $B$ is a function space  with seminorm $\|\cdot\|_B$, we denote by $L^2(0,2\pi;B)$ the class of functions
$u:(0,2\pi)\rightarrow B$ such that
$$
\|u\|_{L^2(B)}:= \left(\Int{0}{2\pi}\|u(t)\|_B^2 {\rm d}t\right)^{\frac 12}<\infty
$$
and similarly, we set
$$
W^{1,2}(0,2\pi;B)=\Big\{u\in L^{2}(0,2\pi;B): \partial_t u\in L^{2}(0,2\pi;B)\Big\}\,.
$$
We will typically shortcut those notations, writing $L^2(B)$ instead of $L^2(0,2\pi;B)$, etc. Next, we define the Banach spaces
\begin{itemize}
\item $L^{q}_\sharp:=\{\bfxi\in L^{q}(0,2\pi) : \bfxi$ is $2\pi$-periodic with $\bar{\bfxi}=\0$\}\,,\ for $q\in[1,\infty]$ with its associated norm $$\|\bfxi\|_{L^q_\sharp}:=\|\bfxi\|_{L^q(0,2\pi)}\,;$$
\item $W^{k}_\sharp:=\{\bfxi\in L^{2}_\sharp(0,2\pi) : d^l{\bfxi}/dt^l\in L^2(0,2\pi)\,,\ \ l=1,\ldots,k\}$ with its associated norm $$\|\bfxi\|_{W^k_\sharp}:=\|\bfxi\|_{W^{k,2}(0,2\pi)}\,;$$
\item $\mathcal L_\sharp^{2}:=\{\bfw\in L^{2}(L^2(\Omega)) : \bfw$ is $2\pi$-periodic, with $\bar{\bfw}=\0\}$ with its associated norm $$\|\bfw\|_{\mathcal L^{2}_\sharp}:=\|\bfw\|_{L^2(L^2(\Omega))}\,;$$
\item $\mathcal W_\sharp^{2}:=\{\bfw\in W^{1,2}(L^2(\Omega))\cap L^2(W^{2,2}(\Omega)) : \bfw$  is $2\pi$-periodic, with $\bar{\bfw}=\0\}$ with its associated norm
$$\|\bfw\|_{\mathcal W_\sharp^{2}}:=\|\bfw\|_{W^{1,2}(L^2(\Omega))}+\|\bfw\|_{L^2(W^{2,2}(\Omega))}\,.$$
\end{itemize}
The function spaces
$$\ba{rl}\medskip
\textsf{{W}$_\sharp^{2}$}&\!\!\!\!:=\left\{\bfw\in L^2(Z^{2,2})\cap W^{1,2}(\calh):\, \mbox{$\bfw$ is $2\pi$-periodic,\,with\, $\bar{\bfw}|_{\Omega_0}=\bar{\hat{\bfw}}=\0$}\right\}\,,\\
\textsf{{L}$_\sharp^{2}$}&\!\!\!\!:=\left\{\bfw\in L^2(\calh):\, \mbox{$\bfw$ is $2\pi$-periodic,\,with\, $\bar{\bfw}|_{\Omega_0}=\bar{\hat{\bfw}}=\0$}\right\}
\ea
$$
are Banach spaces with their corresponding norms
$$
\begin{array}{l}
\|\bfw\|_{\mbox{\scriptsize $\W$}}:=\|\partial_t\bfw\|_{L^2(\Omega)}+\|\bfw\|_{L^2(W^{2,2}(\Omega))}+\|\hat{\bfw}\|_{W_\sharp^1}\,,\\
\|\bfw\|_{\mbox{\scriptsize $\M$}}:=\|\bfw\|_{L^2(L^2(\Omega))}+\|\hat{\bfw}\|_{L^2_\sharp}
\,.
\end{array}
$$
Finally, the natural space for the oscillatory contribution to the pressure is
$$
{\sf P}^{1,2}_\sharp:=\left\{{\sf p}\in L^2(D^{1,2})\ \mbox{with \ $\bar{\sf p}=0$}\right\}\,,
$$
with associated norm
$$
\|{\sf p}\|_{{\sf P}^{1,2}_\sharp}:=\|{\sf p}\|_{ L^2(D^{1,2})}\,.
$$

\section{A General Approach to Time-Periodic Bifurcation}\label{sec:bif}

Consider an abstract evolution problem described by the equation
\be
v_t=M(v,\lambda)\,,
\eeq{Ar.0}
where $\lambda$ is a positive real parameter and $M$ is a nonlinear operator. Suppose that in a neighborhood of $\lambda_o$ $(>0)$, say $I(\lambda_o)$, there exists a branch of steady-state solutions $v_0=v_0(\lambda)$ to \eqref{Ar.0}, namely, for all $\lambda\in I(\lambda_o)$, $v_0(\lambda)$ satisfies
$$
M(v_0,\lambda)=0\,.
$$
The time-periodic bifurcation problem consists then in finding sufficient conditions on $(v(\lambda_o),\lambda_o)$ guaranteeing the existence of a family of time-periodic solutions around $(v_0,\lambda_o)$ that converges to $(v_0,\lambda_o)$ as $\lambda\to\lambda_o$.\par
This problem can be equivalently reformulated in the following, more familiar way. We set
$$u=u(\lambda)=v-v_0(\lambda),\  \mu=\mu(\lambda)=\lambda-\lambda_o\in I(0).$$
It follows from \eqref{Ar.0} that $(u,\mu)$ must solve
\be
u_t+L(u)=N(u,\mu)\,,
\eeq{Ar.1}
where $L$ is the linearization of $M$ at $(v_0,\lambda_o)$  and $N$ is a nonlinear operator depending on the parameter $\mu\in I(0)$, such that $N(0,\mu)=D_uN(0,\mu)=0$ for all such $\mu$. Therefore, the bifurcation problem  at $(v_0(\lambda_o),\lambda_o)$, corresponding to $(u,\mu)=(0,0)$, is equivalent to show that, in addition to the trivial solution $(u,\mu)=(0,0)$, there exists a family of non-trivial $T$--periodic solutions, say $u=u(\mu;t)$, to \eqref{Ar.1} in a neighborhood of $\mu=0$, where the period $T=T(\mu)$ is part of the unkowns, and such that $(u(\mu;t),\mu)\to (0,0)$ as $\mu\to 0$.\par
Setting $\tau:=2\pi\,t/T\equiv \zeta\, t$, \eqref{Ar.1} becomes
\be
\zeta\,u_\tau+L(u)=N(u,\mu)
\eeq{Ar.2}
and the problem reduces to find a family of $2\pi$-periodic solutions to \eqref{Ar.2} with the above properties. If we  write $u=\bar{u}+(u-\bar{u}):=v+w$, where the bar denotes again the average over a $2\pi$-period,  it follows that \eqref{Ar.2} is formally equivalent to the following couple system of equations
\be\ba{ll}\medskip
L(v)=\bar{N(v+w,\mu)}:=N_1(v,w,\mu)\,,\\ \zeta\,w_\tau+L(w)=N(v+w,\mu)-\bar{N(v+w,\mu)}:=N_2(v,w,\mu)\,.\ea
\eeq{Ar.3}

We are thus led to solve the nonlinear ``elliptic-parabolic" system \eqref{Ar.3}. It turns out that, in many circumstances,  the elliptic ``steady-state" problem is better investigated in function spaces with quite less regularity\footnote{Here `regularity' is meant in the sense of behavior at large spatial distances.} than the space where the parabolic ``oscillatory" problem should be studied. This especially occurs  when the physical system evolves in an {\it unbounded spatial region}, in which case the natural framework for the study of \eqref{Ar.3}$_1$ is a {\it homogeneous} Sobolev space, whereas that of \eqref{Ar.3}$_2$ is a classical Sobolev space \cite{GaArma,GaMaH}.
This suggests that, in general, it is more appropriate to study the two equations in \eqref{Ar.3} in two {\it distinct} classes of functions. As a consequence, even though being {\it formally}  the same operator, the operator $L$ in \eqref{Ar.3}$_1$ acts on and ranges into spaces different than those where the operator $L$ in \eqref{Ar.3}$_2$ acts and ranges.\par
With this in mind, \eqref{Ar.3} becomes
\be
L_1(v)=N_1(v,w,\mu)\,,\qquad \zeta\,w_\tau+L_2(w)=N_2(v,w,\mu)\,.
\eeq{Ar.3.2}
The general theory proposed in \cite{GaBif},  and that we are about to recall in its main aspects, takes its cue exactly from these  considerations.\par

In the sequel, we will denote by ${\sf D}[M]\subseteq X$ and ${\sf R}[M]\subseteq Y$, the domain and the  range of a map $M$ between two Banach spaces $X$ and $Y$. We write ${\sf N}[M]:=\{x\in X: M(x)=0\}$ to denote its null space.\par
For a real Banach space $B$, we denote by $B_\comps: = B + {\rm i} B$ its complexification. If $L$ is
a linear operator on $B$, we continue to denote by $L$ its extension to a linear operator on
$X_\comps$, while ${\sf N}_\comps[L]$ and
${\sf R}_\comps[L]$ stand for the null space and range in $B_\comps$. The spectrum, $\sigma(L)$, is computed
with respect to $B_\comps$, so $\nu\in\sigma(L)$ if and only if its complex conjugate, $\nu^*$, is in $\sigma(L)$.\par
Let ${\calu }, {\calv}$,  be real Banach spaces with norms $\|\cdot\|_{\mathcal U}$, $\|\cdot\|_{\mathcal V}$, respectively, and let  $\mathcal Z$ be a real Hilbert space with norm $\|\cdot\|_{\mathcal Z}$.
Moreover,  let
$$
L_1:\mathcal U\mapsto \mathcal V\,
$$
be a bounded linear operator, and let
$$
L_2:{\sf D}\,[L_2]\subset \mathcal Z\mapsto \mathcal Z\,,
$$
be a densely defined, closed linear operator, with a non-empty resolvent set ${\sf P}(L_2)$. For a fixed (once and for all) $\theta\in {\sf P}(L_2)$ we denote by $\mathcal W$ the linear subspace of $\mathcal Z$ closed under the (graph) norm $\|w\|_{\mathcal W}:=\|(L_2+\theta\,I)w\|_{\mathcal Z}$, where $I$ stands for the identity operator in $\calz$.\par
We also define the following spaces
$$\ba{rl}\medskip
\mathcal Z_{2\pi,0}&\!\!\!:=\left\{w:\real\to\calz,\, \mbox{$2\pi$-periodic with} \ \bar{w}=0, \,\mbox{and}\,\int_{0}^{2\pi}\|w(s)\|_{\calz}^2{\rm d}s<\infty \right\}
\\\medskip\calw_{2\pi,0}&\!\!\!:=\left\{w\in \calz_{2\pi,0}:\,\|w\|_{\calw_{2\pi,0}}:=\left(\int_{0}^{2\pi}\left(\|w_t(s)\|_{\calz}^2+\|w(s)\|_{\calw}^2\right){\rm d}s\right)^{\frac12}<\infty\right\}
\ea$$
and we consider a (nonlinear) map
$$
N: \calu\times \calw_{2\pi,0}\times \real\mapsto \calv\oplus \calz_{2\pi,0}
$$
satisfying the following properties:
\be\ba{rl}\medskip
N_1&\!\!\!: (v,w,\mu)\in\calu\times \calw_{2\pi,0}\times \real\mapsto \bar{N(v,w,\mu)}\in \calv
\\
N_2&\!\!\!:=N-N_1:\calu\times \calw_{2\pi,0}\times \real\mapsto \calz_{2\pi,0}\,.
\ea
\eeq{Ar.4}
The Bifurcation Problem is then formulated as follows.\smallskip\par
{\it Find a neighborhood $U(0,0,0)\subset \calu\times \calw_{2\pi,0}\times \real$ of the origin such that the equations}
\be
\ba{ll}\medskip
L_1(v)=N_1(v,w,\mu)\ \mbox{in $\calv$}\,,\\ \zeta\, w_\tau +L_2(w)=N_2(v,w,\mu)\ \mbox{in $\calz_{2\pi,0}$}\,,
\ea
\eeq{Ar.5}
{\it possess there a family of non-trivial $2\pi$-periodic solutions $(v(\mu),w(\mu;\tau))$ for some $\zeta=\zeta(\mu)>0$,  such that $(v(\mu),w(\mu;\cdot))\to (0,0)$ in $\calu\times\calw_{2\pi,0}$ as $\mu\to0$.}\smallskip\par
Whenever the Bifurcation Problem is solvable, we call  $(\mu=0,u=0)$  a {\it bifurcation point}. Moreover, the bifurcation is called {\it supercritical} [resp. {\it subcritical}] if the family of solutions $(v(\mu),w(\mu;\tau))$ exists only for $\mu>0$ [resp. $\mu<0$]. It is called {\it transcritical} otherwise.\par
In order to provide sufficient conditions for the resolution of the above Bifurcation Problem, we introduce the following  assumptions (H1)--(H4) on the involved operators.
\begin{itemize}
  \item[(H1)] $L_1$ is a homeomorphism\,;
  \item[(H2)] There exists $\nu_0:={\rm i}\,\zeta_0$, $\zeta_0>0$ such that
\begin{itemize}
\item[(i)] $L_2-\nu_0I$ is Fredholm of index 0, ${\rm dim}\,{\sf N}_\comps[L_2-\nu_0I]=1$ with ${\sf N}_\comps[L_2-\nu_0I]\cap {\sf R}_{\comps}[L_2-\nu_0I]=\{0\}$. Namely, $\nu_0$ is a simple  eigenvalue of $L_2$.
\item[(ii)]  For every $k\in\nat\backslash\{0,1\}$,  $k\,\nu_0$ is in the resolvent set ${\sf P}(L_2)$ of $L_2$.

\end{itemize}
  \item[(H3)] The operator $Q:\calw_{2\pi,0}\to \calz_{2\pi,0}$ defined by
$$
Q:w\mapsto \zeta_0\,w_\tau+L_2(w)\,,
$$
is Fredholm of index 0\,;
\item[(H4)]
\begin{itemize}
\item[(i)] The operators $N_1,N_2$
are analytic in the neighborhood $U_1(0,0,0)\subset \calu\times \calw_{2\pi,0}\times \real$, i.e.  there exists $\delta>0$ such that for all $(v,w,\mu)\in U_1$ with $\|v\|_{\calu}+\|w\|_{\calw_{2\pi,0}}+|\mu|<\delta$, the Taylor series
$$\ba{ll}\medskip
N_1(v,w,\mu)=\Sum{k,l,m=0}{\infty}R_{klm}v^kw^l\mu^m\,,\\
N_2(v,w,\mu)=\Sum{k,l,m=0}{\infty}S_{klm}v^kw^l\mu^m\,,
\ea
$$
are absolutely convergent in $\calv$ and $\calz_{2\pi,0}$, respectively.
\item[(ii)] The multi-linear operators $R_{klm}$ and $S_{klm}$ satisfy
$$R_{klm}=S_{klm}=0,\quad \text{ whenever }k+l+m\le1,$$
$$R_{011}=R_{00m}=S_{00m}=0,\quad  \text{ for all } m\ge2.$$
\end{itemize}
\end{itemize}

We next define the operator
\be
L_2(\mu):=L_2+\mu\,S_{011}\,,
\eeq{L2mu}
and observe that assumption (i) in (H2) implies that
$$\nu_0\mbox{ is a simple eigenvalue of }L_2(0)\equiv L_2.$$
Under this assumption, denoting the eigenvalues of $L_2(\mu)$ by $\nu(\mu)$, it follows (see e.g.  \cite[Proposition 79.15 and Corollary 79.16]{Z1}) that the map $\mu\mapsto\nu(\mu)$ is well defined and of class $C^\infty$ in a neighborhood of $\mu=0$.
\renewcommand{\theequation}{\arabic{section}.\arabic{equation}}
\par
\setcounter{equation}{8}
In the sequel, for all $z\in\C$, we set:
$$\Re[z]\mbox{ its real part, }\ \Im[z]\mbox{  its imaginary part, }\ {z}^*\mbox{ its conjugate.}
$$
\smallskip\par
We may now state the following bifurcation result which will be the key to formulate \theoref{8.1} for problem \eqref{04}.

\Bt\label{3.1_ar} Suppose  {\rm (H1)--(H4)} hold and, in addition,
\be
\Re[\nu'(0)]\neq 0\,.\footnote{This is the usual transversality condition : the  eigenvalue $\nu(\mu)$ crosses the imaginary axis with ``non-zero speed."}
\eeq{nupr}
 Moreover, let
 \begin{itemize}
 \item $v_0$ be a normalized eigenvector of $L_2$ corresponding to the eigenvalue $\nu_0$,
 \item $v_1:=\Re[v_0\,{\rm e}^{-{\rm i}\,\tau}].$
\end{itemize}
Then, the following properties are valid. \smallskip\\
{\rm (a)} {\rm Existence.} There is a neighborhood $\mathcal I(0)$ of $0\in\real$, there are analytic families
\be
\big(v(\varepsilon),w(\varepsilon),\zeta(\varepsilon),\mu(\varepsilon)\big)\in \calu\times \calw_{2\pi,0}\times \real_+\times\real
\eeq{fam}
satisfying \eqref{Ar.5}, for all $\varepsilon\in\mathcal I(0)$, and such that
\be
\big(v(\varepsilon),w(\varepsilon)-\varepsilon\,v_1,\zeta(\varepsilon),\mu(\varepsilon)\big)\to (0,0,\zeta_0,0)\ \ \mbox{as $\varepsilon\to 0$}\,.
\eeq{Ar.10}
\par\noindent
{\rm (b)} {\rm Uniqueness.}
There is a neighborhood  $$U(0,0,\zeta_0,0)\subset \calu\times \calw_{2\pi,0}\times \real_+\times \real$$ such that every (nontrivial) $2\pi$-periodic solution to \eqref{Ar.5},  lying in $U$ must coincide, up to a phase shift, with a member of the family \eqref{fam}.
\smallskip\par\noindent
{\rm (c)} {\rm Parity.}  The functions $\zeta(\varepsilon)$ and $\mu(\varepsilon)$ are even:
$$
\zeta(\varepsilon)=\zeta(-\varepsilon)\,,\ \ \mu(\varepsilon)=\mu(-\varepsilon)\,,\ \ \mbox{for all $\varepsilon\in\cali(0)$\,.}
$$\label{3.1_ar}
\ET{3.1_ar}

As a consequence of {\rm (c)}, the bifurcation due to these solutions is either subcritical or supercritical, a two-sided bifurcation being excluded, unless $\mu\equiv 0$. \theoref{3.1_ar} was shown in \cite{GaBif} but, for completeness and reader's sake, we deem it appropriate to present a proof in the Appendix.

\section{Analysis of the linearized operators}\label{chap:bifu}

The existence
of a branch of equilibria of \eqref{02} in a class of homogeneous Sobolev spaces, for arbitrary values of the Reynolds number $\lambda\, (>0)$ follows easily from classical results regarding steady-state Navier-Stokes problems in exterior domains. Indeed, we know from \cite[Theorem X.6.4]{Gab} that for any $\lambda>0$  problem \eqref{03}$_{1-4}$ has  one corresponding  solution $(\bfu_0,p_0)$ in the class \eqref{sfaco} recalled below. We then set
\be
\bfchi_0:= -\varpi\mathbb{A} ^{-1}\int_{\partial\Omega}\mathbb T(\bfu_0,p_0)\cdot\bfn\,,
\eeq{2.3_00}
which is well defined by standard trace theorems. This yields

\Bt For any $\lambda>0$, problem \eqref{03} has at least one solution
$${\sf s}_0(\lambda):=(\bfu_0(\lambda),p_0(\lambda),\bfchi_0(\lambda))$$ such that
\be {\sf s}_0(\lambda)\in [L^q(\Omega)\cap D^{1,r}(\Omega)\cap D^{2,s}(\Omega)]\times[L^{\sigma}(\Omega)\cap D^{1,s}(\Omega)]\times\real^3,
\eeq{sfaco}
for all $q\in (2,\infty]$, $r\in(\frac43,\infty]$, $s\in (\frac32,\infty]$, $\sigma\in (1,\infty)$.
\ET{exi}
More properties on equilibrium configurations, such as uniqueness and stability, are studied in \cite{BoGaGa1}. In a nutshell, it is proved therein that for small values of $\lambda$, the equilibrium is unique, locally stable and an oscillatory motion does not exist.

\subsection{Spectral Properties}\label{sec:spectrum}
We look, from now on, for solutions $(\bfu,p,\bfchi)$ to the problem \eqref{02}, that we recall here for convenience
\be\ba{cc}\medskip\left.\ba{r}\medskip
\partial_t\bfu-\Delta\bfu+\nabla p=\lambda\,[\partial_1\bfu+(\dot{\bfchi}-\bfu)\cdot\nabla\bfu]\\
\Div\bfu=0\ea\right\}\ \ \mbox{in $\Omega\times(0,\infty)$}\,,\\ \medskip
\bfu(x,t)={\dot{\bfchi}}(t)+\bfe_1\,, \ \mbox{ $(x,t)\in\partial\Omega\times(0,\infty)$}\,;\ \
\Lim{|x|\to\infty}\bfu(x,t)=\0\,,\ t\in(0,\infty)\,,\\
\ddot{\bfchi}+\mathbb{A} \cdot \bfchi+\varpi\Int{\partial\Omega}{} \mathbb T(\bfu,p)\cdot\bfn=\0 \ \ \mbox{in $(0,\infty)$}\,,
\ea
\eeq{02-chap4}
bifurcating from a branch of equilibria ${\sf s}_0(\lambda)=(\bfu_0(\lambda),p_0(\lambda),\bfchi_0(\lambda))$ of \eqref{02-chap4}, i.e. a branch of solutions of
\be\ba{cc}\medskip\left.\ba{r}\medskip
-\Delta\bfu_0+\nabla p_0=\lambda\,(\partial_1\bfu_0-\bfu_0\cdot\nabla\bfu_0)\\
\Div\bfu_0=0\ea\right\}\ \ \mbox{in $\Omega$}\,,\\ \medskip
\bfu_0(x)=\bfe_1\,, \ \mbox{ $x\in\partial\Omega$}\,;\ \
\Lim{|x|\to\infty}\bfu_0(x)=\0\,,\\
\mathbb{A} \cdot \bfchi_0+\varpi\Int{\partial\Omega}{} \mathbb T(\bfu_0,p_0)\cdot\bfn=\0\,.
\ea
\eeq{03-chap4}
Writing
$$
\bfu=\bfu_0(\lambda)+\bfw(t;\lambda)\,,\ p=p_0(\lambda)+{\sfp}(t;\lambda)\,,\ \ \bfchi=\bfchi_0(\lambda)+ \bfeta(t;\lambda)\,,
$$
equations \eqref{02-chap4} yield the first-order in time system with respect to the four unknowns $(\bfw,{\sfp},\bfeta,\bfsigma)$:
\be\ba{cc}\medskip\left.\ba{lr}\medskip
& \partial_t\bfw+\lambda\,[-\partial_1\bfw+\bfu_0\Cdot\nabla\bfw+
({\bfsigma}-\bfw)\Cdot\nabla\bfu_0 +({\bfsigma}-\bfw)\Cdot\nabla\bfw] =\Delta\bfw-\nabla {\sfp}\\ \medskip
& \Div\bfw=0\ea\right\}\,  \mbox{in $\Omega\times(0,\infty)$}\,,\\ \medskip
\bfw(x,t)={{\bfsigma}}(t)\,, \ \mbox{ $(x,t)\in\partial\Omega\times(0,\infty)$}\,;\ \
\Lim{|x|\to\infty}\bfw(x,t)=\0\,,\ t\in(0,\infty)\,,\\
\left.\ba{r}\medskip
\dot{\bfsigma}+\mathbb{A} \cdot \bfeta+\varpi\Int{\partial\Omega}{} \mathbb T(\bfw,{\sfp})\Cdot\bfn=\0\,\\
\dot{\bfeta}(t)=\bfsigma(t)\ea\right\}\,  \mbox{in $(0,\infty)$}\,. &
\ea
\eeq{04}

We aim to reformulate this bifurcation problem in a suitable functional setting, which will eventually allow us to apply the theory presented in Section~\ref{sec:bif}. The major  challenge in reaching this goal is the choice of the {\it right setting}, where the involved operators possess the properties required by \theoref{3.1_ar}. This will require a careful study of appropriate linearized problems obtained from \eqref{03-chap4} and \eqref{04} that will be performed in the present and following subsections.\par
Here we secure some relevant spectral properties of the operator obtained by linearizing \eqref{04} around the trivial solution
$$(\bfw,{\sfp},\bfsigma,\bfeta)\equiv(\0,0,\0,\0).$$
In this regard, we begin to define the maps
$$
\tilde{\Delta}:(\bfw,\bfeta)\in Z^{2,2}\times\real^3\mapsto \tilde{\Delta}(\bfw,\bfeta)\in \call^2(\real^3)\,,
$$
where $Z^{2,2}$ was introduced in \eqref{zetaa},
\be
\tilde{\Delta}(\bfw,\bfeta)=\left\{\ba{ll}\medskip
 -\Delta\bfw\ \mbox{in $\Omega$}\,,\\
\mathbb{A} \cdot \bfeta+2\varpi\Int{\partial\Omega}{}\mathbb D(\bfw)\cdot\bfn\ \ \mbox{in $\Omega_0$\,,}\ea\right.
\eeq{Delta}
and
$$
\tilde{\partial_1}:\bfw\in Z^{2,2}\mapsto \tilde{\partial_1}(\bfw)\in \calh(\real^3)\,,
$$
where
\be
\tilde{\partial_1}(\bfw)=\left\{\ba{ll}\medskip
 -\partial_1\bfw\ \mbox{in $\Omega$}\,,\\
\0\ \ \mbox{in $\Omega_0$\,.}\ea\right.
\eeq{D1}

One can readily check that $\tilde{\partial}_1(\bfw)\in \calh(\real^3)$. Set
\be
\mathscr A:=\mathscr P\,\tilde{\Delta},\quad  {\bsfW}:=(\bfw,\bfeta),\quad \bfsigma:=\bfw|_{\Omega_0},
\eeq{defAbsWbfsigma}
where $\mathscr P$ is  the orthogonal projection of $\call^2(\real^3)$ onto $\calh(\real^3)$, see e.g. \cite[Lemma 2]{BoGaGa1}. Then we consider the operator defined by
\arraycolsep=1pt
\be
\ba{rl}
\mathscr L_0:  Z^{2,2}\times \real^3\subset \calh(\real^3)\times\real^3&\to\calh(\real^3)\times\real^3\\
{\bsfW} & \mapsto(\lambda_o\tilde{\partial}_1\bfw+\mathscr A(\bfw,\bfeta),-\bfsigma).
\ea
\eeq{L0}
\arraycolsep=5pt
\par
In the sequel, $\cali$ denotes the identity in $\calh_{\mbox{\tiny $\comp$}}(\real^3)\times\mathbb C^3$.
The following result holds. \Bl Let $\zeta\in\real\backslash\{0\}$. Then, the operator $\mathscr L_0-\i\,\zeta\cali$ is a homeomorphism from $Z_{\mbox{\tiny $\mathbb C$}}^{2,2}\times\comp^3$ onto $\calh_{\mbox{\tiny $\comp$}}(\mathbb R^3)\times\comp^3$.  Moreover,
if $|\zeta|\ge1$,
there exists $c=c(\Omega_0,\mathbb{A},\varpi,\lambda_o)>0$, such that
\be
\|D^2\bfw\|_2+|\zeta|^{\frac12}\|\nabla\bfw\|_2+|\zeta|(\|\bfw\|_2+|\bfsigma|+|\bfeta|)\le c\,\|(\mathscr L_0-\i\,\zeta\cali)(\bsfW)\|_2\,.
\eeq{4.2}
\EL{4.1}
\begin{proof} By the orthogonal decomposition given in \cite[Lemma 2]{BoGaGa1},  \eqref{Delta} and \eqref{D1}, we infer that, given $\bfcalf:=((\bff,\bfF),\bfG)\in\calh_{\mbox{\tiny $\mathbb C$}}(\real^3)\times\comp^3$, the equation $$(\mathscr L_0-\i\,\zeta\cali)(\bsfW)=\bfcalf$$
is equivalent to the following
set of equations
\be\ba{cc}\medskip\left.\ba{r}\medskip\Delta\bfw+\lambda_o\partial_1\bfw-\nabla{\mathfrak p}=-\i\,\zeta\,\bfw -\bff\\
\Div\bfw=0\ea\right\}\,  \mbox{in $\Omega$}\,,\\ \medskip
\bfw=\bfsigma\ \ \mbox{on $\partial\Omega$}\,,\\ \medskip
-\i\,\zeta\,{\bfsigma}+\mathbb{A}\cdot \bfeta+\varpi\Int{\partial\Omega}{}\mathbb T(\bfw,{\mathfrak p})\cdot\bfn=\bfF\,,\\ \medskip
\bfsigma=-\i\,\zeta\bfeta-\bfG\,.
\ea
\eeq{4.1}
Therefore, the homeomorphism property follows if, for any given $((\bff,\bfF),\bfG)$ specified above, problem \eqref{4.1}
has one and only one solution $(\bfw,\mathfrak p,\bfsigma,\bfeta)\in W_{\mbox{\tiny $\mathbb C$}}^{2,2}(\Omega)\times D^{1,2}_{\mbox{\tiny $\mathbb C$}}(\Omega)\times\mathbb C^3\times\mathbb C^3$. We begin to establish some formal estimates, holding for all $\zeta\neq 0$. In this regards, we recall that, by \cite[Lemmas 9--11]{Gah},
\be
|\bfsigma|\le c_0\|\mathbb D(\bfw)\|_2
\eeq{tra}
with $c_0=c_0(\Omega_0)>0$.
Dot-multiplying both sides of \eqref{4.1}$_1$ by $\bfw^*$ and integrating by parts over $\Omega$, taking account of \eqref{4.1}$_{2-5}$, we infer
\be
\begin{split}
2\|\mathbb D(\bfw)\|_2^2-\i\,\zeta\,\left(\|\bfw\|_2^2+\frac1\varpi|\bfsigma|^2-\frac{1}{\varpi}\bfeta^*\cdot\mathbb{A}\cdot\bfeta \right)-\lambda_o(\partial_1\bfw,\bfw^*)\\ =(\bff,\bfw^*)+\frac1\varpi\bfF\cdot\bfsigma^*+\frac{1}{\varpi}\bfG^*\cdot\mathbb{A}\cdot\bfeta\,.
\end{split}
\eeq{4.3}
Taking the real part of \eqref{4.3},  using Korn identity, Cauchy--Schwarz inequality, and the fact that $\Re \,(\partial_1\bfw,\bfw^*)=~0$, we deduce
\be
\|\nabla\bfw\|_2^2\le c_1\left(\|\bff\|_2\|\bfw\|_2+|\bfF||\bfsigma|+|\bfG||\bfeta|\right),
\eeq{4.4}
where, here and in the rest of the proof, $c_1$ denotes a positive constant depending, at most, on $\Omega_0$, $\mathbb{A}$, $\lambda_o$, and $\varpi$.\par
From \eqref{4.1}$_5$ and \eqref{tra} we deduce
\be
|\bfeta|\le {c_1}{|\zeta|^{-1}}\left(|\bfG|+\|\nabla\bfw\|_2\right),
\eeq{xi}
so that employing in \eqref{4.4} the latter, \eqref{tra}, Korn identity and Cauchy-Schwarz inequality we show
$$
\|\nabla\bfw\|_2^2\le c_1\left[(|\zeta|^{-1}+|\zeta|^{-2})|\bfG|^2+\|\bff\|_2\|\bfw\|_2+|\bfF||\bfsigma|\right].
$$
If we apply Cauchy inequality on the last two terms on the right-hand side of this inequality, we may deduce, on the one hand,
\be
\|\nabla\bfw\|_2^2\le c_1\left[|\zeta|^{-1}(|\bfF|^2+|\bfG|^2+\|\bff\|_2^2)+|\zeta|^{-2}|\bfG|^2\right]+|\zeta|(|\bfsigma|^2+\|\bfw\|_2^2)
\eeq{Du0}
and, on the other hand,
\be
\|\nabla\bfw\|_2^2\le c_1\left[|\zeta|^{-1}|\bfG|^2+|\zeta|^{-2}\left(|\bfF|^2+|\bfG|^2+\|\bff\|_2^2\right)\right]+\varepsilon |\zeta|^2(|\bfsigma|^2+\|\bfw\|_2^2)\,.
\eeq{Du}
where $\varepsilon>0$ is arbitrarily fixed and $c_1$ depends also on $\varepsilon$.\par
We now take the imaginary part of \eqref{4.3} and use Cauchy--Schwarz inequality to get
$$
|\zeta|(|\bfsigma|^2+\|\bfw\|_2^2)\le c_1\left[|\zeta||\bfeta|^2+\|\nabla\bfw\|_2\|\bfw\|_2+\|\bff\|_2\|\bfw\|_2+|\bfG||\bfeta|+|\bfF||\bfsigma|\right]\,.
$$
By utilizing, on the right-hand side of the latter, \eqref{xi} and Cauchy inequality we obtain
\be
|\zeta|(|\bfsigma|^2+\|\bfw\|_2^2)\le {c_1}\,|\zeta|^{-1}\left[\|\nabla\bfw\|_2^2+|\bfG|^2+|\bfF|^2+\|\bff\|_2^2\right]\,.
\eeq{ima}
Thus, combining \eqref{Du} and \eqref{ima}, and taking $\varepsilon$ suitably small we deduce
\be
|\zeta|^2(|\bfsigma|^2+\|\bfw\|_2^2)\le {c_1}(1+|\zeta|^{-1}+|\zeta|^{-2})\left(|\bfG|^2+|\bfF|^2+\|\bff\|_2^2\right)\,,
\eeq{ima1}
which, in turn, once replaced in \eqref{Du0}, delivers
\be
|\zeta|\|\nabla\bfw\|_2^2\le {c_1}(1+|\zeta|^{-1}+|\zeta|^{-2})\left(|\bfG|^2+|\bfF|^2+\|\bff\|_2^2\right)\,.
\eeq{Du1}
Finally, from \eqref{xi} and \eqref{Du1} we conclude
\be
|\zeta|^2|\bfeta|^2\le c_1(1+|\zeta|^{-1}+|\zeta|^{-2}+|\zeta|^{-3})\left(|\bfG|^2+|\bfF|^2+\|\bff\|_2^2\right)\,.
\eeq{xi1}

By combining estimates \eqref{ima1}--\eqref{xi1} with the classical Galerkin method, one obtains that for any given $((\bff,\bfF),\bfG)$ in the specified class and $\zeta\neq 0$, there exists a (unique, weak) solution to \eqref{4.1} such that
$$(\bfw,\bfeta)\in [\cald^{1,2}_{\mbox{\tiny $\mathbb C$}}(\real^3)\cap\call^2_{\mbox{\tiny $\mathbb C$}}(\real^3)]\times\compl^3,$$
satisfying \eqref{ima1}--\eqref{xi1}. We next write \eqref{4.1}$_{1-3}$ in the following Stokes-problem form
$$\ba{cc}\medskip\left.\ba{r}\medskip\Delta\bfw=\nabla{\mathfrak p}+\bfcalg\\
\Div\bfw=0\ea\right\}\,  \mbox{in $\Omega$}\,,\\ \medskip
\bfw=\bfsigma\ \ \mbox{on $\partial\Omega$}
\ea
$$
where
$$
\bfcalg:=-\lambda_o\, \partial_1\bfw +\i\,\zeta\,\bfw-\bff\,.
$$
Since $\bfcalg\in L^2_{\mbox{\tiny $\mathbb C$}}(\Omega)$ and $\bfu\in W^{1,2}_{\mbox{\tiny $\mathbb C$}}(\Omega)$, from classical results \cite[Theorems IV.5.1 and V.5.3]{Gab} it follows that $D^2\bfw\in L^2(\Omega)$, thus completing the existence (and uniqueness) proof.\par Furthermore, by \cite[Lemmas IV.1.1 and V.4.3]{Gab} we get
\be
\|D^2\bfw\|_2\le c\,\Big[\|\bff\|_2+(\lambda_o+1)\|\nabla\bfw\|_2+(|\zeta|+1)\,\|\bfw\|_2+|\bfsigma|\Big]\,.
\eeq{4.7}
As a result, if $|\zeta|\ge 1$, the inequality in \eqref{4.2} is a consequence of \eqref{ima1}--\eqref{4.7}, and the proof of the lemma is completed.
\end{proof}\par
Let $\bsfu_0\in X^{2}(\Omega)$ and consider the operator
\be
\hat{\mathscr K}: \bsfW:=(\bfw,\bfeta)\in  Z^{2,2}\times\real^3 \mapsto \hat{\mathscr K}(\bsfW)\in \call^2(\real^3)
\eeq{3.32}
where, recalling also \eqref{defAbsWbfsigma},
\be
\hat{\mathscr K}(\bsfW)=\left\{\ba{ll}\medskip
\lambda_o\big(\bsfu_0\cdot\nabla\bfw+(\bfw-\bfsigma)\cdot\nabla\bsfu_0\big)\ \ \mbox{in $\Omega$}\,,\\
\0\ \ \mbox{in $\Omega_0$}\,.
\ea\right.
\eeq{3.32_1}
Recall that the space $X(\Omega)$ is continuously embedded in $L^4(\Omega)$, see e.g. \cite[Proposition 5]{BoGaGa1}.
Then, classical results imply
\be\bsfu_0\in L^4(\Omega)\cap D^{1,2}(\Omega)\,,\ \ W^{2,2}(\Omega)\subset  L^\infty(\Omega)\cap W^{1,4}(\Omega)\,,
\eeq{u0w}
so that the operator $\mathscr K$ is well defined.
Next, let
\be
\mathscr K: \bsfW:=(\bfw,\bfeta)\in  Z^{2,2}\times\real^3 \mapsto (\mathscr P\hat{\mathscr K},\0)\in \calh(\real^3)\times\real^3\,,
\eeq{DeFi}
and define
$\mathscr L_2:Z^{2,2}\times \real^3\to \calh(\real^3)\times\real^3$  by
\be
\mathscr L_2(\bsfW):=\mathscr L_0(\bsfW) +\mathscr K(\bsfW).
\eeq{L2}
The main result of this subsection reads as follows.
\Bt For all $\zeta\neq0$, the operator
$
\mathscr L_2-\i\,\zeta\cali$ is Fredholm of index 0. Moreover, denoting by $\sigma(\mathscr L_2)$  the spectrum of $\mathscr L_2$, we have that
\begin{itemize}
\item $\sigma(\mathscr L_2)\cap \{\i\,\real\backslash\{0\}\}$
consists  of an at most finite or countable number of eigenvalues, each of which is isolated and of finite (algebraic) multiplicity;
\item $\sigma(\mathscr L_2)\cap \{\i\,\real\backslash\{0\}\}$
has $0$ as only cluster point.
\end{itemize}
\ET{5.1}
\begin{proof} We begin to prove  that the operator $\hat{\mathscr K}$ defined in \eqref{3.32}--\eqref{3.32_1} is compact. Let $(\bsfW_k)_k$ be a sequence bounded in $Z^{2,2}\times\real^3$. This implies, in particular, that there is $M>0$  independent of $k$ such that
\be
\|\bfw_k\|_{2,2}+|\bfsigma_k|\le M\,.
\eeq{5.9}
Then, the compact embedding $W^{2,2}(\Omega)\subset W^{1,4}(\Omega_R)\cap L^\infty(\Omega_R)$, for all $R>1$,  ensures the existence of $(\bfw_*,\bfsigma_*)\in W^{2,2}(\Omega)\times\real^3$ and subsequences, again denoted by $(\bfw_k,\bfsigma_k)_k$, such that
\be
\begin{split}
\bfw_k\to\bfw_*\ & \mbox{strongly in $W^{1,4}(\Omega_R)\cap L^\infty(\Omega_R)$, for all $R>1$}\,;\\ \bfsigma_k\to\bfsigma_*\ & \mbox{in $\real^3$}\,.
\end{split}
\eeq{JN}
In view of the linearity of $\hat{\mathscr K}$, we may assume, without loss of generality,  $\bfw_*\equiv\bfsigma_*\equiv\0$. By H\"older inequality, we deduce
$$\ba{ll}\medskip
\|\hat{\mathscr K}(\bsfW_k)\|_2\le \lambda_o\left[\|\bsfu_0\|_4\|\bfw_k\|_{1,4,\Omega_R}+\|\bsfu_0\|_{4,\Omega^R}\|\bfw_k\|_{1,4}\right.\\
\hspace*{3.5cm}+\left.\|\nabla\bsfu_0\|_2(\|\bfw_k\|_{\infty,\Omega_R}+|\bfsigma_k|)+\|\nabla\bsfu_0\|_{2,\Omega^R}\|\bfw_k\|_{\infty}\right]\,.
\ea
$$
Therefore, letting $k\to\infty$ in this inequality, and using \eqref{u0w},  \eqref{5.9}, and  \eqref{JN}, we infer, with a positive constant $C$ independent of $R$, that
$$
\lim_{k\to\infty}\|\hat{\mathscr K}(\bsfW_k)\|_2\le C\,\left(\|\bsfu_0\|_{4,\Omega^R}+\|\nabla\bsfu_0\|_{2,\Omega^R}\right)\,.
$$
This, in turn, again by \eqref{u0w} and the arbitrariness of $R>1$ shows the desired property. The compactness of $\hat{\mathscr K}$ and \lemmref{4.1} then imply that the operator
$$
\widehat{\mathscr L}_\zeta:=\mathscr L_2-\i\,\zeta\cali
$$
is Fredholm of index 0, for all $\zeta\neq0$. The theorem is then a consequence of well-known results (e.g. \cite[Theorem XVII.4.3]{GG}) provided we show that the null space of $\widehat{\mathscr L}_\zeta$ is trivial, for all sufficiently large $|\zeta|$. To this end, we observe that $\widehat{\mathscr L}_\zeta(\bsfW)=0$ means $\mathscr L_0(\bsfW)-\i\,\zeta\,\bsfW=-\mathscr K(\bsfW)$. Applying \eqref{4.2} and recalling \eqref{DeFi}, we thus get
$$
|\zeta|^{\frac12}\|\nabla\bfw\|_2+|\zeta|(\,\|\bfw\|_2+|\bfsigma|)\le c\,\|\hat{\mathscr K}(\bsfW)\|_2\,,\ \text{ for all }  |\zeta|\ge 1,
$$
where $c$ is independent of $\zeta$. Also, from \eqref{3.32_1} and H\"older inequality, we deduce
$$
\|\hat{\mathscr K}(\bsfW)\|_2\le \lambda_o\, \|\bsfu_0\|_{1,\infty}(\|\bfw\|_{1,2}+|\bfsigma|)\,.
$$
combining the last two inequalities, and since \theoref{exi} implies $\|\bsfu_0\|_{1,\infty}<\infty$, we conclude that $\bfw\equiv\0$ when $|\zeta|$ is larger than a suitable positive constant. This completes the proof of the theorem.
\end{proof}\par

\subsection{The Linearized Time-Periodic Operator}\label{sec:time_per}
Objective of this subsection is to show that Assumption ({H3}) of \theoref{3.1_ar} holds in our setting. This will be the content of \theoref{3.1}.

To this end, we begin to demonstrate some fundamental functional properties of the linear operator obtained by formally setting $(\bfsigma-\bfw)\cdot\nabla\bfw\equiv\0$ in equations \eqref{04}, in the class of periodic solutions with zero average. To do this,  we need to investigate the more special linear problem to which \eqref{04} reduces by taking $(\bfsigma-\bfw)\cdot\nabla\bfw\equiv\bfu_0\equiv\0$, that is,
\be\ba{cc}\medskip\left.\ba{lr}\medskip
& \partial_t\bfw-\lambda\partial_1\bfw=\Delta\bfw-\nabla {\sfp} \\
& \Div\bfw=0\ea\right\}\ \ \mbox{in $\Omega\times(0,\infty)$}\,,\\ \medskip
\bfw(x,t)={{\bfsigma}}(t)\,, \ \mbox{ $(x,t)\in\partial\Omega\times(0,\infty)$}\,;\ \
\Lim{|x|\to\infty}\bfw(x,t)=\0\,,\ t\in(0,\infty)\,,\\
\left.\ba{r}\medskip
\dot{\bfsigma}+\mathbb{A} \cdot \bfeta+\varpi\Int{\partial\Omega}{} \mathbb T(\bfw,{\sfp})\Cdot\bfn=\0\,\\
\dot{\bfeta}(t)=\bfsigma(t)\ea\right\}\,  \mbox{in $(0,\infty)$}\,. &
\ea
\eeq{04-simple}
\par
To reach our goal,
we begin to prove some crucial  properties of
the associated family of boundary-value problems for Fourier modes. Precisely, for  $k\in\mathbb Z$, and $m\in\{1,2,3\}$, we introduce  the pair of fields $(\bfh^{(m)}_k,p_k^{(m)})$, such that $(\bfh_0^{(m)},p^{(m)}_0)\equiv (\0,0)$ and, for $k\neq0$,
\be
\ba{cc}\medskip\left.\ba{r}\medskip
{\rm i}\,k\,\zeta_0\,\bfh_k^{(m)}-\lambda_o\partial_1\bfh_k^{(m)}=\Delta \bfh_k^{(m)}-\nabla p_k^{(m)}\\
\Div \bfh_k^{(m)}=0\ea\right\}\ \ \mbox{in $\Omega$}\\
\bfh_k^{(m)}=\zeta_0\bfe_m\ \ \mbox{on $\partial\Omega$}\,,
\ea
\eeq{3.2}
with $\zeta_0\in\real\setminus\{0\}$ and $\lambda_o\in\real$. Moreover, for a fixed $k$, we denote by $\mathbb K$ the matrix whose components are given by
\be
({\mathbb K})_{\ell m}=\left(\IdS\mathbb T(\bfh^{(m)}_k,p^{(m)}_k)\cdot\bfn\right)_\ell\,,\quad \ell,m=1,2,3.
\eeq{Matrices}
\par Unless otherwise stated, throughout the rest of this subsection  we shall
use summation convention over repeated indices.\smallskip
\par
The following result plays a key role in our analysis.

\Bl Let $\zeta_0\in\real\setminus\{0\}$, $\lambda_o>0$.\footnote{In fact, the lemma continues to hold for any $\lambda_o\in\real$, but this is irrelevant to our aims.} For every $m\in\{1,2,3\}$ and every $k\in \mathbb Z\backslash\{0\}$, \eqref{3.2} has a unique solution $(\bfh_k^{(m)},p_k^{(m)})\in W^{2,2}(\Omega)\times D^{1,2}(\Omega)$ which, moreover, satisfies the estimates
\be
\left\|\bfh_k^{(m)}\right\|_2\le C\,\left\|\nabla\bfh_k^{(m)}\right\|_2\le C \,(|k|+1)^{\frac12}\,,\qquad
\left\|D^2\bfh_k^{(m)}\right\|_2\le C \,(|k|+1)\,,
\eeq{3.4}
where $C$ is a constant independent of $k$. In addition, the following properties hold :
\begin{itemize}
\item[{\rm (i)}]  $\mathbb K$  is invertible and $\mathbb{A}-k^2\,\zeta_0^2\mathbb I +{\rm i}\,k\,\varpi\,\mathbb K$  is invertible for every $k\in \mathbb Z$ ;
\item[{\rm (ii)}] for every $\bfalpha\in \mathbb C^3$, we have
\be
{\rm i}\,k\,\zeta_0\,\|\bfsf h_k\|_2^2+2\|\mathbb D(\bfsf h_k)\|_2^2-\lambda_o(\partial_1\bfsf h_k,\bfsf h_k^*)=\zeta^2_0\bfalpha^*\cdot\mathbb K\cdot\bfalpha
\eeq{3.A}
where $\bfsf h_k:= \alpha_m\bfh^{(m)}_k$.
\end{itemize}
\EL{3.1}
\begin{proof}
Since the proof is the same
for $m = 1, 2, 3$, we pick $m = 1$ and will omit the superscript.
Let
\be
\bsfe(x):=\zeta_0\curl(x_2\phi(|x|)\bfe_3)
\eeq{ee}
where $\phi=\phi(|x|)$ is a smooth ``cut-off" function that is 1 in $\Omega_{\rho_1}$ and 0 in $\bar{\Omega^{\rho_{2}}}$, $R_*<\rho_1<\rho_2$. Clearly, $\bsfe$ is smooth in $\Omega$, with bounded support and, in addition, $\bsfe(x)=\zeta_0\bfe_1$ in a neighborhood of $\partial\Omega$ and $\Div\bsfe=0$ in $\Omega$.
From \eqref{3.2} we then deduce that $\bfv_k:=\bfh_k-\bsfe$ solves the following boundary-value problem  (for all $|k|\ge1$) :
\be
\ba{cc}\medskip\left.\ba{r}\medskip
{\rm i}\,k\,\zeta_0\,\bfv_k-\lambda_o\,\partial_1\bfv_k=\Delta \bfv_k-\nabla p_k+\lambda_o\partial_1\bsfe-{\rm i}\,k\,\zeta_0\,\bsfe+\Delta\bsfe\\
\Div \bfv_k=0\ea\right\}\ \ \mbox{in $\Omega$}\\
\bfv_k=\0\ \ \mbox{on $\partial\Omega$}\,.
\ea
\eeq{3.5}

We shall next show a number of a priori estimates for solutions $(\bfv_k,p_k)$  to \eqref{3.5} in their stated function class that once combined, for instance, with the classical Galerkin method, will produce the desired existence  result. Uniqueness will also be  an obvious consequence of these estimates.
For this, we dot-multiply both sides of \eqref{3.5}$_1$ by $\bfv^*_k$, and integrate by parts to get
\be
{\rm i}\,k\,\zeta_0\,\|\bfv_k\|_2^2-\lambda_o(\partial_1\bfv_k,\bfv_k^*)+\|\nabla\bfv_k\|_2^2=(\bfcalf_k,\bfv_k^*)\,,
\eeq{3.6}
where $\bfcalf_k:=\lambda_o\partial_1\bsfe-{\rm i}\,k\,\zeta_0\bsfe+\Delta\bsfe$. By the properties of $\bsfe$,
\be
\|\bfcalf_k\|_2\le c\,(|k|+1)
\eeq{3.7}
where, here and in the rest of the proof, $c$ denotes a generic (positive) constant depending, at most, on $\zeta_0,\lambda_o$ and $\Omega$, but {\it not} on $k$. Since
\be
\Re\, (\partial_1\bfv_k,\bfv_k^*)=0\,,
\eeq{3.8}
by taking the real part of \eqref{3.6} and using \eqref{3.7} we infer
\be
\|\nabla\bfv_k\|_2^2\le c\,(|k|+1)\|\bfv_k\|_2\,.
\eeq{3.9}

Considering the imaginary part of \eqref{3.6} and employing \eqref{3.7}--\eqref{3.9} along with Cauchy--Schwarz inequality, we obtain
$$
|k|\|\bfv_k\|_2\le c\,(\|\nabla\bfv_k\|_2+|k|+1)\le c(|k|+1)\|\bfv_k\|_2^{\frac12}\,,
$$
which gives the uniform bound, independent of $k$,
\be
\|\bfv_k\|_2\le c\,.
\eeq{3.10}
Taking into account that $\bfh_k=\bsfe+\bfv_k$, \eqref{3.10} proves \eqref{3.4}$_1$. Inequality \eqref{3.4}$_2$ then follows by using  \eqref{3.10} into \eqref{3.9}.\par
Moreover, from classical estimates on the Stokes problem, see e.g.\ \cite[Lemma 1]{Hey}, we infer that
$$
\|D^2\bfv_k\|_2\le c\,\big(\|\partial_1\bfv_k\|_2+\|\bfcalf_k\|_2+\|\nabla\bfv_k\|_2\big)
$$
and so, recalling that $|k|\ge 1$,  from the latter  inequality, \eqref{3.7}, \eqref{3.9} and \eqref{3.10} we show \eqref{3.4}$_3$.
\par
Let $\bfalpha\in\mathbb C^3$, and, for fixed $k\neq0$,
 set
$$
\bfsf h_k:= \alpha_m\bfh^{(m)}_k\,,\ \ {\sfp}_k:=\alpha_m\,p^{(m)}_k\,.
$$
From \eqref{3.2}  we then find
\be
\ba{cc}\medskip\left.\ba{r}\medskip
{\rm i}\,k\,\zeta_0\bfsf h_k-\lambda_o\partial_1\bfsf h_k=\Div\mathbb T(\bfsf h_k,{\sfp}_k)\\
\Div \bfsf h_k=0\ea\right\}\ \ \mbox{in $\Omega$}\\
\bfsf h_k=\zeta_0\bfalpha\,\ \mbox{on $\partial\Omega$}.
\ea
\eeq{3.2_ar}
Dot-multiplying both sides of \eqref{3.2_ar}$_1$ by $\bfsf h^*_k$ and integrating by parts over $\Omega$ we deduce
\be
\i\,k\,\zeta_0\|\bfsf h_k\|_2^2+\|\mathbb D(\bfsf h_k)\|_2^2-\lambda_o(\partial_1\bfsf h_k,\bfsf h^*_k)=\zeta^2_0\bfalpha^*\cdot\mathbb K\cdot\bfalpha\,,
\eeq{previousrelation}
that is, \eqref{3.A}.\par
Now, suppose that there exists $\hat{\bfalpha}\in\mathbb C^3$ such that
\[\mathbb{A}\cdot\hat{\bfalpha}-k^2\,\zeta_0^2\hat{\bfalpha} +{\rm i}\,k\,\varpi\,\mathbb K \cdot\hat{\bfalpha} =0.\]
This implies that either $k=0$ in which case $\hat{\bfalpha}=0$ because $\mathbb{A}$ is invertible, or
$$\hat{\bfalpha}^*\cdot\mathbb K \cdot\hat{\bfalpha}= -{\rm i}\,\hat{\bfalpha}^*\cdot\mathbb{L}\cdot\hat{\bfalpha} \, ,$$
for some real valued matrix $\mathbb{L}$. Then from \eqref{previousrelation}, we infer that
$$
\i\left(k\,\zeta_0\|\hat{\bfsf h}_k\|_2^2+\hat{\bfalpha}^*\cdot\mathbb{L}\cdot\hat{\bfalpha} \right)-\lambda_o(\partial_1\hat{\bfsf h_k},\hat{\bfsf h}^*_k)=2\|\mathbb D(\hat{\bfsf h_k})\|_2^2\,,
$$
which, in turn, recalling that \be
\Re\,(\partial_1\hat{\bfsf h_k},\hat{\bfsf h}^*_k)=0\,,
\eeq{re0}
allows us to we deduce $\hat{\bfsf h}_k=\0$ in $W^{2,2}(\Omega)$. The latter implies $\hat{\bfalpha}=\0$ and thus shows the desired property for $\mathbb K$ when $\mathbb L=0$ and for $\mathbb{A}-k^2\,\zeta_0^2\mathbb I +{\rm i}\,k\,\varpi\,\mathbb K$ otherwise.  The proof of the lemma is completed.
\end{proof}\par
With the help of \lemmref{3.1}, we deduce the following statement.

\Bc Let $\zeta_0\in\real\setminus\{0\}$, $\lambda_o>0$, and $\varpi\in \real\backslash\{0\}$. Then, for any $\bsfF\in  L^2_\sharp$, the  problem
\be\ba{cc}\medskip\left.\ba{r}\medskip\zeta_0\partial_\tau\bsfw-\lambda_o\,\partial_1\bsfw=\Delta\bsfw-\nabla{\sf q}\\
\Div\bsfw=0\ea\right\}\,  \mbox{in $\Omega\times [0,2\pi]$}\,,\\ \medskip
\bsfw=\zeta_0\dot{\bfeta}\ \ \mbox{on $\partial\Omega\times [0,2\pi]$}\,,\\
\zeta_0^2\ddot{\bfxi}+\mathbb{A}\cdot\bfxi+\varpi\Int{\partial\Omega}{}\mathbb T(\bsfw,{\sf q})\cdot\bfn=\bsfF\,,\, \ \mbox{in $[0,2\pi]$}\,,
\ea
\eeq{3.1}
has one and only one solution $\big(\bsfw,{\sf q},\bfxi\big)\in \mathcal W_\sharp^{2}
\times {\sf P}^{1,2}_\sharp   \times W^{2}_\sharp$. Moreover, we have the estimate
\be
\|\bsfw\|_{\mathcal W_\sharp^{2}}+\|{\sf q}\|_{{\sf P}^{1,2}_\sharp}+\|\bfxi\|_{W^{2}_\sharp}\le C\,\|\bsfF\|_{L^2_\sharp}\,,
\eeq{n}
for some $C=C(\Omega,\lambda_o,\mathbb{A},\zeta_0)>0$.
\EC{3.2_0}
\begin{proof}
We formally expand $\bsfw$, ${\sf q}$, and  $\bfxi$ in Fourier series as follows:
\be
\bsfw(x,t)=\Sum{k\in\mathbb Z}{}\bsfw_k(x)\,{\rm e}^{\i k\,t}\,,\qquad{\sf q}(x,t)=\Sum{k\in\mathbb Z}{}{\sf q}_k(x)\,{\rm e}^{\i k\,t}\,,\qquad
\bfxi(t)=\Sum{k\in\mathbb Z}{}\bfxi_k \,{\rm e}^{\i k\,t}\,,
\eeq{Fou}
where $\bsfw_0\equiv\nabla{\sf q}_0\equiv\bfxi_0\equiv\0$ and for $k\neq0$,
 $(\bsfw_k,{\sf q}_k,\bfxi_k)$ solve the problem
\be\ba{cc}\medskip\left.\ba{r}\medskip
\i\,k\,\zeta_0\,\bsfw_k-\lambda_o\partial_1\bsfw_k=\Delta \bsfw_k-\nabla{\sf q}_k\\
\Div\bsfw_k=0\ea\right\}\ \ \mbox{in $\Omega$}\\
\bsfw_k|_{\partial\Omega}=\zeta_0\,{\rm i} k \bfxi_k\,,
\ea
\eeq{2.11}
with the further condition
\be
\left(-k^2\,\zeta_0^2\mathbb I+\mathbb{A} \right)\cdot \bfxi_k+\varpi\Int{\partial\Omega}{}\mathbb T(\bsfw_k,{\sf q}_k)\cdot\bfn=\bsfF_k\,,
\eeq{2.8_1}
where $(\bsfF_k)_k$ are Fourier coefficients of $\bsfF$ (so that also $\bsfF_0\equiv\0$).
For each fixed $k\in\mathbb Z\setminus\{0\}$, a  solution to \eqref{2.11}--\eqref{2.8_1} is given by\footnote{No summation over $k$ here!}
\be
\bsfw_k=\sum_{\ell=1}^3\i\,k\,\bfxi_{k\ell}\bfh_k^{(m)}\,,\ \ {\sf q}_k=\sum_{\ell=1}^3\i\,k\,\bfxi_{k\ell}p_k^{(m)}\,,
\eeq{1.20}
with $(\bfh_k^{(m)},p_k^{(m)})$ given in \lemmref{3.1}, and where $\bfxi_k$ solve the equations
\be
\left(-k^2\,\zeta_0^2\mathbb I+\mathbb{A}\right)\cdot \bfxi_k+\sum_{\ell=1}^3{\rm i}\,k\,\varpi\,\bfxi_{k\ell}\left(\Int{\partial\Omega}{}\mathbb T(\bfh_k^{(m)},{p}^{(m)}_k)\cdot\bfn\right)_\ell=\bsfF_k\,,
\eeq{chi}
which, with the notation \eqref{Matrices}, can be equivalently rewritten as
\be
\mathbb M\cdot\bfxi_k=\bsfF_k\,,\ \ \mathbb M:=(-k^2\,\zeta_0^2\mathbb I+\mathbb{A} )+{\rm i}\,k\,\varpi\,\mathbb K\,.
\eeq{1.21}
By \lemmref{3.1}--(i), the matrix $\mathbb M$ is invertible for all $k\neq 0$.
As a result, for the given $\bsfF_k$, equation  \eqref{1.21} has one and only one solution $\bfxi_k$.
If we now dot-multiply both sides of \eqref{1.21} by $\bfxi_k^*$ and use  \eqref{3.A}, we deduce
$$\ba{ll}\medskip
-k^2\,\zeta_0^2|\bfxi_k|^2+\bfxi_k^*\cdot\mathbb{A} \cdot\bfxi_k\,-\zeta^{-2}_0\Big[k^2\zeta_0\,\varpi\,
\|\bfxi_{km}\bfh_k^{(m)}\|_2^2\\ \hspace*{1cm}+ {\rm i}\lambda_o\varpi\,(\partial_1(\bfxi_{km}\bfh_k^{(m)}),(\bfxi_{km}\bfh_k^{(m)})^*)
-2{\rm i}\,k\,\varpi\,\|\mathbb D(\bfxi_{km}\bfh_k^{(m)})\|_2^2\Big]=(\bsfF_k,\bfxi_k^*)\,,
\ea
$$
which, in view of \eqref{re0}, furnishes
\be\ba{rl}\medskip
2k\,\zeta_0^{-2}\,\varpi\|\mathbb D(\bfxi_{km}\bfh_k^{(m)})\|_2^2&\!\!\!\!=\Im\,(\bsfF_k,\bfxi_k^*)\
\ea
\eeq{2.17-1}
and
\be\ba{rl}\medskip
-k^2\,\zeta_0^2|\bfxi_k|^2+\bfxi_k^*\cdot\mathbb{A} \cdot\bfxi_k\,-\zeta_0^{-2}\Big[k^2\zeta_0\varpi\,
\|\bfxi_{km}\bfh_k^{(m)}\|_2^2+\\ {\rm i}\lambda_o\varpi\,(\partial_1(\bfxi_{km}\bfh_k^{(m)}),(\bfxi_{km}\bfh_k^{(m)})^*)\Big]&\!\!\!\!=\Re\,(\bsfF_k,\bfxi_k^*)\,.
\ea
\eeq{2.17-2}
Recalling that $\bfxi_{km}\bfh_k^{(m)}|_{\partial\Omega}=\bfxi_k$, by \eqref{2.17-1}, Cauchy--Schwarz inequality and trace inequality, we obtain
the crucial estimate
\be
|\bfxi_k|+\|\nabla(\bfxi_{km}\bfh_k^{(m)})\|_2\le c\,|k|^{-1}\,|\bsfF_k|\,,\ \ |k|\ge1\,,
\eeq{crucial}
where, here and in the following, $c$ denotes a generic positive constant independent of $k$.
Again by Cauchy--Schwarz inequality, from \eqref{2.17-2} we get
$$
k^2\,(|\bfxi_k|^2+\|\bfxi_{km}\bfh_k^{(m)}\|_2^2)\le c\,\left(|\bfxi_k|^2+\|\bfxi_{km}\bfh_k^{(m)}\|_2\|\nabla(\bfxi_{km}\bfh_k^{(m)})\|_2+|\bsfF_k||\bfxi_k|\right)
$$
from which, using \eqref{crucial} and Cauchy inequality
we deduce
$$
k^2\,(|\bfxi_k|^2+\|\bfxi_{km}\bfh_k^{(m)}\|_2^2)\le c\,|k|^{-2}|\bsfF_k|^2\,,\ \ |k|\ge1\,,
$$
that allows us to conclude
$$
k^4|\bfxi_k|\le c\,|\bsfF_k|^2\,, \ \ |k|\ge1.
$$
It immediately follows that
\be
\|\bfxi\|^2_{W^{2}_\sharp}=\sum_{|k|\ge1}(|k|^4+|k|^2+1)|\bfxi_k|^2\le c\sum_{|k|\ge1}|\bsfF_k|^2= c\|\bsfF\|_{L^2_\sharp}^2\,.
\eeq{1.46}
Moreover, we  infer from \eqref{1.20},
\eqref{1.46} and \eqref{3.4} that
\be
\begin{split}
\|\bsfw\|_{\calw_\sharp^2}^2 & =\!\sum_{|k|\ge1}\left[(|k|^2+1)\|\bsfw_k\|_2^2+\|\nabla\bsfw_k\|_2^2+\|D^2\bsfw_k\|_2^2\right]\\ & \le c\!\sum_{|k|\ge1}(|k|^4+|k|^2+1)|\bfsigma_k|^2\le c\,\|\bsfF\|_{L^2_\sharp}^2\,,
\end{split}
\eeq{1.49}
so that, combining \eqref{1.46},  \eqref{1.49}, and  \eqref{3.1}$_{1}$, we obtain
\be
\|\bsfw\|_{\calw_\sharp^2}+\|\bfxi\|_{W^{2}_\sharp}+\|{\sf q}\|_{{\sf P}^{1,2}_\sharp}\le c\,\|{\bsfF}\|_{L^2_\sharp}\,.
\eeq{2.22}
and the proof of existence is completed.\par
The uniqueness property amounts to show that the problem
\be\ba{cc}\medskip\left.\ba{r}\medskip
\zeta_0\partial_\tau\bsfw-\lambda_o\,\partial_1\bsfw=\Delta\bsfw-\nabla {\sf q}\\
\Div\bsfw=0\ea\right\}\ \ \mbox{in $\Omega\times[0,2\pi]$}\\ \medskip
\bsfw|_{\partial\Omega}=\dot{\bfeta}\,;\\\zeta_0^2 \ddot{\bfeta}+\mathbb{A} \cdot \bfeta+\varpi\Int{\partial\Omega}{}\mathbb T(\bsfw,{\sf q})\cdot\bfn=\0\
\ea
\eeq{3.31}
has only the zero solution in the specified function class. If we dot-multiply \eqref{3.31}$_1$ by $\bsfw$, integrate by parts over $\Omega$ and use \eqref{3.31}$_3$, we get
\[
\half\ode{}t(\zeta_0\|\bsfw(t)\|_2^2+\zeta_0^2|\dot{\bfeta}(t)|^2+\bfeta(t)\cdot \mathbb{A}\cdot \bfeta(t)|^2)+2\|\mathbb D(\bsfw(t))\|_2^2=0\,.
\]
Integrating both sides of this equation from $0$ to $2\pi$ and employing the $2\pi$-periodicity of the solution, we easily obtain  $\|\mathbb D(\bsfw(t))\|_2\equiv 0$ which in turn immediately furnishes\footnote{Recall that $\bsfw(t)\in L^2(\Omega)$.}
$\bsfw\equiv\nabla{\sf q}\equiv\0$. The proof of the corollary is completed.\end{proof}

\Br The invertibility of the matrix $\mathbb M$, defined in \eqref{1.21}, is crucial to the resolution of \eqref{chi} and, therefore, to the result stated in the corollary above. The remarkable fact is that this property holds  for {\it all} values of $\zeta_0>0$ and $|k|\ge1$, independently of the eigenvalues of the matrix $\mathbb{A}$. In physical terms, this means that the possibility of a ``disruptive" resonance is always ruled out. However, observe that  if there is some $\bar{k}$ such that $\bar{k}^2\zeta_0=\omega_{\sf n}^2$, where $\omega_{\sf n}^2$ is an eigenvalue of $\mathbb  A$, then the invertibility of $\mathbb M$ obviously degenerates in the limit $\varpi\to 0$. To simplify the discussion, assume for instance that $\mathbb{A}=\omega_{\sf n}^2\mathbb I$. Then it follows from \eqref{1.21} that the corresponding amplitude $|\bfxi_{\bar{k}}|$ satisfies
$$
|\bfxi_{\bar{k}}|=\frac{\sqrt{\zeta_0}}{\omega_{\sf n}\varpi}|\mathbb K^{-1}\cdot\bsfF_{\bar{k}}|\,,
$$
which may become increasingly large when $\varpi\to 0$, namely, when the mass of the liquid volume occupied by $\mathscr B$  is vanishingly small compared to that of $\mathscr B$.
\ER{9.1}
\par
The next lemma proves well-posedness of the linear problem \eqref{04-simple}, in the class of $2\pi$-periodic solutions with zero average. The crucial aspect of this result is that there is no need to impose restrictions on $\zeta_0$ with respect to the natural frequencies of the restoring force.
\Bl Let $\lambda_o,\zeta_0, \varpi \in (0,\infty)$.\footnote{We could, more generally,  assume $\zeta_0\in\real\backslash\{0\}$ which, again, would be immaterial; see \cite[Lemma 1.1]{You}.} Then, for any $(\bff,\bfF,\bfG)\in \call_\sharp^{2}\times L^2_\sharp\times W_\sharp^1$, the  problem
\be\ba{cc}\medskip\left.\ba{r}\medskip\zeta_0\partial_\tau\bfw-\lambda_o\partial_1\bfw=\Delta\bfw-\nabla{p}+\bff\\
\Div\bfw=0\ea\right\}\,  \mbox{in $\Omega\times [0,2\pi]$}\,,\\ \medskip
\bfw=\zeta_0\,\dot{\bfxi}-\bfG\ \ \mbox{on $\partial\Omega\times [0,2\pi]$}\,,\\
\,\zeta_0^2\ddot{\bfxi}+\mathbb{A}  \cdot \bfxi+\varpi\Int{\partial\Omega}{}\mathbb T(\bfw,{p})\cdot\bfn=\bfF\,, \ \ \mbox{in $[0,2\pi]$}\,,
\ea
\eeq{3.1}
has one and only one solution $\big(\bfw,p,\bfxi\big)\in \mathcal W_\sharp^{2}
\times {\sf P}^{1,2}_\sharp\times W^{2,2}_\sharp$. This solution satisfies the estimate
\be
\|\bfw\|_{\mathcal W_\sharp^{2}}+\|{p}\|_{{\sf P}^{1,2}_\sharp}+\|\bfxi\|_{W^{2,2}}\le C\,\Big(\|\bff\|_{\call_\sharp^{2}}+\|\bfF\|_{L^2_\sharp}+\|\bfG\|_{W^1_\sharp}\Big)\,,
\eeq{n}
where $C=C(\Omega,\lambda_o,\zeta_0,\mathbb{A},\varpi)$.
\EL{3.2}
\begin{proof} Let $\bfw=\bfz+\bfu$ where $\bfz$ and $\bfu$ satisfy the following problems
\be\ba{cc}\medskip\left.\ba{r}\medskip
\zeta_0\partial_\tau\bfz-\lambda_o\partial_1\bfz=\Delta\bfz-\nabla {\sf r}+\bff\\
\Div\bfz=0\ea\right\}\ \ \mbox{in $\Omega\times [0,2\pi]$}\\
\bfz|_{\partial\Omega}=-\bfG
\ea
\eeq{2.3}
and
\be\ba{cc}\medskip\left.\ba{r}\medskip\zeta_0^2
\partial_\tau\bfu-\lambda_o\partial_1\bfu=\Delta\bfu-\nabla {\sf q}\\
\Div\bfu=0\ea\right\}\ \ \mbox{in $\Omega\times[0,2\pi]$}\\ \medskip
\bfu|_{\partial\Omega}=\zeta_0\,\dot{\bfxi}\,;\\ \medskip\zeta_0^2\ddot{\bfxi}+\mathbb{A} \cdot \bfxi+\varpi\Int{\partial\Omega}{}\bfT(\bfu,{\sf q})\cdot\bfn=\bfF-\varpi\Int{\partial\Omega}{}\bfT(\bfz,{\sf r})\cdot\bfn:=\bsfF\,,\ \ \mbox{in $[0,2\pi]$}\,.
\ea
\eeq{3.17}
Set
\be
\bfW(x,t):=x_3G_2(t)\bfe_1+x_1G_3(t)\bfe_2+x_2G_1(t)\bfe_3\,,
\eeq{W}
so that
\be
\curl\bfW=\bfG(t)\,.
\eeq{W1}
Let $\phi(x)$ be the ``cut-off" function defined in the beginning of the proof of \lemmref{3.1}, and define
$$
\bfw(x,t):=\curl\big(\phi(x)\bfW(x,t)\big)\,.
$$
In view of \eqref{W1} we deduce
\be
\bfw(x,t)=\phi(x)\bfG(t)-\bfW(x,t)\times\nabla\phi(x)
\eeq{03_1}
so that $\bfw$ is a $2\pi$-periodic solenoidal vector function
that is equal to $\bfG(t)$  for $|x|\le \rho_1$ and equal to $0$ for $|x|\ge \rho_2$. Therefore, from \eqref{2.3} we deduce that the field
\be
\bsfz(x,t):=\bfz(x,t)+\bfw(x,t)\,,
\eeq{zeta}
obeys the following problem
\be\ba{cc}\medskip\left.\ba{r}\medskip
\zeta_0\partial_\tau\bsfz-\lambda_o\partial_1\bsfz=\Delta\bsfz-\nabla {\sf r}+\bsff\\
\Div\bsfz=0\ea\right\}\ \ \mbox{in $\Omega\times [0,2\pi]$}\\
\bsfz|_{\partial\Omega}=\0\,,
\ea
\eeq{2.3_1}
where
\be
\bsff:=\bff-\zeta_0\partial_\tau\bfw+\lambda_o\partial_1\bfw-\Delta\bfw\,.
\eeq{EDB}
From \eqref{W}--\eqref{03_1},  \eqref{EDB} and the assumption on $\bfG$ it follows that $\bsff\in \call^2_\sharp$ and that
\be
\|\bsff\|_{\call_\sharp^2}\le c\,(\|\bff\|_{\call_\sharp^2}+\|\bfG\|_{W^2_\sharp})\,.
\eeq{mds}

Employing \cite[Theorem 12]{GaMaH}, we then deduce  that there exists a unique solution $(\bsfz,{\sf r})\in  \mathcal W_\sharp^{2}\times{\sf P}^{1}_\sharp$   that, in addition, obeys the inequality
$$
\|\bsfz\|_{\mathcal W_\sharp^{2}}+\|{\sf r}\|_{{\sf P}^{1,2}_\sharp}\le c\,\|\bsff\|_{\mathcal L_\sharp^{2}}\,.
$$
The latter, in combination with \eqref{mds} and \eqref{zeta}, allows us to conclude $\bfz\in\calw_\sharp^2$ and
\be
\|\bfz\|_{\mathcal W_\sharp^{2}}+\|{\sf r}\|_{{{\sf P}^{1,2}_\sharp}}\le c\,(\|\bff\|_{\mathcal L_\sharp^{2}}+\|\bfG\|_{W^2_\sharp})\,.
\eeq{3.18}
Now, by the trace theorem\footnote{Possibly, by modifying $\sf r$ by adding to it a suitable function of time.} and \eqref{3.18} we get
$$
\left\|\Int{\partial\Omega}{}\mathbb T(\bfz,{\sf \sf r})\cdot\bfn\right\|_{L^2}\le c\,\left(\|\bfz\|_{\mathcal W^{2}_\sharp}+\|{\sf \sf r}\|_{{\sf P}^{1,2}_\sharp}\right)\le c\,(\|\bff\|_{\mathcal L_\sharp^{2}}+\|\bfG\|_{W^2_\sharp})\,,
$$
so that  $\bsfF$ belongs to $L^2_\sharp(0,2\pi)$ and satisfies
\be
\|\bsfF\|_{L^2_\sharp}\le c(\|\bff\|_{\call^2_\sharp}+\|\bfF\|_{L^2_\sharp}+\|\bfG\|_{W^2_\sharp})
\,.
\eeq{3.19}
Thus, from \cororef{3.2_0} it follows that there is one and only one solution $(\bfu,{\sf q},\bfxi)\in \calw_\sharp^2\times {{\sf P}^{1,2}_\sharp}\times W^{2,2}_\sharp$ to \eqref{3.17} that, in addition, satisfies the estimate
$$
\|\bfu\|_{\mathcal W_\sharp^{2}}+\|{\sf q}\|_{{\sf P}^{1,2}_\sharp}+\|\bfxi\|_{W^{2}_\sharp}\le c\,\|\bsfF\|_{L^2_\sharp}\,.
$$
As a result, combining the latter with \eqref{3.19} and \eqref{3.18}, we
 complete the existence proof.

The stated uniqueness property amounts to show that every solution $(\bfw,{\sfp},\bfxi)$ to \eqref{3.31} in the stated function class vanishes identically and this has been already shown at the end of the proof of \cororef{3.2_0}.
\end{proof}\par
As said at the beginning of this subsection, our ultimate aim is to  show that, within our functional setting,  Assumption ({H3}) is satisfied.
For this purpose we define yet another suitable linearized operator that will allow us to
formulate \eqref{3.1}, with $\bfF:=\bff|_{\Omega_0}$, as the operator equation:
\be
\mathscr Q_0(\bsfW)=\bsfF,
\eeq{chiappa}
with $\bsfF:=(\bff,\bfG)\in \M\times W^1_\sharp$.
This amounts to define $\mathscr Q_0:{\W}\times W^2_\sharp\to \M\times W^1_\sharp$ by
$$
\mathscr Q_0:\bsfW\mapsto \mathscr \zeta_0\partial_t\bsfW+\mathscr L_0(\bsfW)\,,
$$
where $\mathscr L_0$ is given by \eqref{L0}, i.e.
\be
\mathscr L_0(\bsfW)= (\lambda_o\tilde{\partial}_1\bfw+\mathscr A(\bfw,\bfxi),-\bfsigma)\,,
\eeq{L0-2}
with $\sigma=\bfw|_{\Omega_0}$. We emphasize that, since
$$\mathscr Q_0(\bsfW)=\zeta_0\partial_t\bsfW+\mathscr L_0(\bsfW) \equiv (\zeta_0\partial_t \bfw + \lambda_o\tilde{\partial}_1\bfw+\mathscr A(\bfw,\bfxi),\zeta_0\dot \bfxi-\bfsigma)\,,
$$
the second component of the equation
$$\mathscr Q_0(\bsfW)=\bsfF,$$  gives
the side condition $\zeta_0\dot \bfxi-\bfw = \bfG$ in $\eqref{3.1}_3$.
Then, from \lemmref{3.2} we infer the following important result.
\Bl For any given $\lambda_o, \zeta_0,\varpi\in (0,\infty)$, the operator
$$
\mathscr Q_0:{\W}\times W^2_\sharp\to \M\times W^1_\sharp
$$
is a homeomorphism.
\EL{Hom}
\par
Next, define $\mathscr Q:{\W}\times W^2_\sharp\to \M\times W^1_\sharp$ by
$$
\mathscr Q:\bsfW\mapsto \zeta_0\partial_t\bsfW+\mathscr L_2(\bsfW)\,,
$$
where $\mathscr L_2$ has the form given in \eqref{L2}, i.e.
\be
\mathscr L_2(\bsfW)=\mathscr L_0(\bsfW) +\mathscr K(\bsfW),
\eeq{L2-2}
where, we recall, $\mathscr K(\bsfW)=(\mathscr P\hat{\mathscr K}(\bsfW),\0)$,  and
\be
\hat{\mathscr K}(\bsfW)=\left\{\ba{ll}\medskip
\lambda_o\big(\bsfu_0\cdot\nabla\bfw+(\bfw-\bfsigma)\cdot\nabla\bsfu_0\big)\ \ \mbox{in $\Omega$}\,,\\
\0\ \ \mbox{in $\Omega_0$}\,,
\ea\right.
\eeq{3.32_2}
with $\bfsigma=\bfw|_{\Omega_0}$.

\lemmref{Hom} allows us to prove the main finding of this subsection, namely,   that Assumption ({H3}) is satisfied. Precisely, we have the following result.
\Bt For any given $\lambda_o, \zeta_0,\varpi\in (0,\infty)$,
the operator $\mathscr Q$ is Fredholm of index 0.
\ET{3.1}
\begin{proof} Since $\mathscr Q=\mathscr Q_0+\mathscr K$,  by \lemmref{Hom}  the stated property will follow provided we show that the map
$$
\hat{\mathscr K} :\bsfW\in\W\mapsto \hat{\mathscr K}(\bsfW)\in \M
$$
defined as in \eqref{3.32_2} is compact. Let $(\bsfW_k)_k$ be a bounded sequence in $\W$. This means  that there is $M>0$ independent of $k$ such that
\be
\|\bfw_k\|_{\calw_\sharp^2}+\|\bfsigma_k\|_{W_\sharp^{1}}\le M\,,
\eeq{gatto1}
with $\bfsigma_k:=\bfw_k|_{\Omega_0}$.
We may then select sequences (still denoted by $(\bfw_k,\bfsigma_k)_k$) such that
\be
\bfw_k\to{\bfw_*} \ \, \mbox{weakly in $\W$\,;}\ \ \bfsigma_k\to{\bfsigma_*}\,  \  \mbox{strongly in $L^{\infty}(0,2\pi)$,}
\eeq{2.18}
for some $(\bfw_*,\bfsigma_*)\in \mathcal W^{2}_\sharp\times W^{1,2}_\sharp$.
Due to the linearity of $\hat{\mathscr K}$, without loss of generality  we may take $\bfw_*\equiv\bfsigma_*\equiv\0$, so that we must  show that
\be
\lim_{k\to\infty}\int_0^{2\pi}\|\mathscr K(\bsfW_k)\|_{2,\Omega}^2=0\,.
\eeq{gatto2}

From \eqref{2.18},  the compact embeddings $W^{2,2}(\Omega)\subset W^{1,4}(\Omega_R)$ for all $R>R_*$, and Lions-Aubin lemma we infer
\be
\int_{0}^{2\pi}\left(\|\bfw_k(t)\|_{4,\Omega_R}^2+\|\nabla\bfw_k(t)\|_{4,\Omega_R}^2\right)\to 0\ \ \mbox{as $k\to\infty$, for all $R>R_*$\,.}
\eeq{2.19}
Furthermore,
$$
\int_{0}^{2\pi} \|\bsfu_0\cdot\nabla\bfw_k(t)\|_{2,\Omega}^2
\le \|\bsfu_0\|_4^2\int_{0}^{2\pi} \|\nabla\bfw_k(t)\|_{4,\Omega_R}^2+ \|\bsfu_0\|_{4,\Omega^R}^2\int_{0}^{2\pi}\|\nabla\bfw_k(t)\|_{4,\Omega}^2\,,
$$
which, by the embeddings $X(\Omega)\subset L^4(\Omega)$ and  $W^{2,2}(\Omega)\subset W^{1,4}(\Omega)$,  \eqref{gatto1}, \eqref{2.19}  and the arbitrariness of $R$ furnishes
\be
\lim_{k\to\infty}\int_{0}^{2\pi} \|\bsfu_0\cdot\nabla\bfw_k(t)\|_{2}^2=0\,.
\eeq{2.21}
Likewise,
$$
\Int{0}{2\pi} \|\bfw_k(t)\cdot\nabla\bsfu_0\|_{2,\Omega}^2
\le \|\nabla\bsfu_0\|_{4}^2\Int{0}{2\pi} \|\bfw_k(t)\|_{4,\Omega_R}^2 +
\|\nabla\bsfu_0\|_{4,\Omega^R}^2\Int{0}{2\pi}\|\bfw_k(t)\|_{4,\Omega}^2\,,
$$
so that, arguing as before,
\be
\lim_{k\to\infty}\int_{0}^{2\pi} \|\bfw_k(t)\cdot\nabla\bsfu_0\|_{2,\Omega}^2=0\,.
\eeq{2.22}

Finally,
$$
\int_{0}^{2\pi}\|\bfsigma_k(t)\cdot\nabla\bsfu_0\|_{2,\Omega}^2\le {2\pi}\, \|\bfsigma_k\|_{L^\infty(0,2\pi)}^2\|\nabla\bsfu_0\|_2^2\,,
$$
so that \eqref{2.18}$_2$ yields
\be
\lim_{k\to\infty}\int_{0}^{2\pi}\|\bfsigma_k(t)\cdot\nabla\bsfu_0\|_{2,\Omega}^2=0
\eeq{3.38}
Combining \eqref{2.21}--\eqref{3.38} we thus arrive at \eqref{gatto2},
which completes the proof of the theorem.\end{proof}

\section{The Bifurcation Result}
This section contains our main results, \theoref{NEC}  and \theoref{8.1}. We start by rewriting the bifurcation problem as a system of operator equations in an appropriate functional setting that, thanks to the results proved in the previous sections,  will enable us to apply \theoref{3.1_ar}.
\subsection{Reformulation of the Problem in Banach Spaces}\label{sec:rif}
The first step is to split \eqref{04} into its averaged and oscillatory components. Set
$$
\begin{array}{lll}
\bfw:=\bar{\bfw}+(\bfw-\bar{\bfw}):=\bsfu+\bsfw\,,\qquad &\bfeta=\bar{\bfeta}+(\bfeta-\bar{\bfeta}):=\bar{\bfeta}+\bfxi\,,\\ {\sfp}=\bar{\sfp}+({\sfp}-\bar{\sfp}):=\bar{\sfp}+{\sf q}\,,\qquad &\mu:=\lambda-\lambda_o\,,\\
\bsfu_0:=\bfu_0(x;\lambda_o)\,,\qquad &\tilde{\bfu}_0:=\bfu_0(x;\mu+\lambda_o)-\bsfu_0.
\end{array}$$
From \eqref{04} we thus get
\be\ba{cc}\medskip\left.\ba{r}\medskip
-\lambda_o(\partial_1\bsfu-\bsfu_0\cdot\nabla\bsfu-\bsfu\cdot\nabla\bsfu_0)-\Delta\bsfu+\nabla\bar{\sfp}=\bfN_1(\bsfu,\bsfw,\mu)\\
\Div\bsfu=0\ea\right\}\ \ \mbox{in $\Omega$\,,}\\ \medskip
\bsfu=\0\ \ \mbox{on $\partial\Omega$}\,,\\
\mathbb{A}  \cdot \bar{\bfeta}+\varpi\Int{\partial\Omega}{}\mathbb T(\bsfu,\bar{\sfp})\cdot\bfn=\0\,,
\ea
\eeq{7.1}
where
\be
\bfN_1:=\left\{\ba{ll}\ms-\mu(\partial_1\bsfu-\bsfu_0\cdot\nabla\bsfu-\bsfu\cdot\nabla\bsfu_0)\\-(\mu+\lambda_o)\left[\tilde{\bfu}_0\cdot\nabla\bsfu+\bsfu\cdot\nabla\tilde{\bfu}_0+\bsfu\cdot\nabla\bsfu+\bar{(\bfsigma-\bsfw)\cdot\nabla\bsfw}\right]\ \mbox{in $\Omega$}\\
\0\ \mbox{in $\Omega_0$}
\,,
\ea\right.
\eeq{7.2}
and, with the time-scaling $\tau:=\zeta\,t$,
\be\ba{cc}\medskip\left.\ba{r}\medskip
\zeta\partial_\tau\bsfw-\lambda_o(\partial_1\bsfw-\bsfu_0\cdot\nabla\bsfw+(\bfsigma-\bsfw)\cdot\nabla\bsfu_0)-\Delta\bsfw+\nabla{\sf q} =\bfN_2(\bsfu,\bsfw,\mu)\\
\Div\bsfw=0\ea\right\}\ \ \mbox{in $\Omega\times [0,2\pi]$}\\ \medskip
\bsfw=\bfsigma\ \ \mbox{on $\partial\Omega\times[0,2\pi]$}\\
\left.\ba{r}\medskip
\zeta\dot{\bfsigma}+\mathbb{A}  \cdot \bfxi+\varpi\Int{\partial\Omega}{}\mathbb T(\bsfw,{\sf q})\cdot\bfn=\0\\ \zeta\dot{\bfxi}-\bfsigma=\0\ea\right\}\,  \mbox{in $(0,\infty)$}\,.
\ea
\eeq{7.3}
where
\be
\bfN_2:=\left\{\ba{ll}\medskip\ba{ll}\medskip\mu(\partial_1\bsfw-\bsfu_0\cdot\nabla\bsfw+(\bfsigma-\bsfw)\cdot\nabla\bsfu_0)\\
+(\mu+\lambda_o)\Big[\tilde{\bfu}_0\cdot\nabla\bsfw+(\bsfw-\bfsigma)\cdot\nabla\tilde{\bfu}_0+\bar{(\bfsigma-\bsfw)\cdot\nabla\bsfw}\\
+{(\bfsigma-\bsfw)\cdot\nabla\bsfw}+\bsfu\cdot\nabla\bsfw+(\bfsigma-\bsfw)\cdot\nabla\bsfu\Big]\ \mbox{in $\Omega$}\ea\\
\0 \ \mbox{in $\Omega_0$}.
\ea\right.
\eeq{7.4}

We first determine the image of $\bfN_1$ and $\bfN_2$.
\Bl
Let $\bfu_0=\bfu_0(\lambda)$ be the velocity field of the solution determined in \theoref{exi},  corresponding to $\lambda>0$ and let $(\bsfu,\bsfw,\mu)\in X^2(\Omega)\times {\sf W}^2_\sharp\times\real$. Then, the following properties hold:
\be
\bfN_1(\bsfu,\bsfw,\mu)\in \mathcal Y(\Omega)\,;\ \ \bfN_2(\bsfu,\bsfw,\mu)\in {\sf L}^2_\sharp\,.
\eeq{7.5}
\EL{7.1}
\begin{proof} Taking into account   \eqref{7.2} and arguing as in the proof of \cite[Lemma 16]{BoGaGa1}, it is easy to verify that the first property in \eqref{7.5} holds if
\be
\bar{(\bfsigma-\bsfw)\cdot\nabla\bsfw}\in  \cald_0^{-1,2}(\Omega)\cap L^2(\Omega)\,.
\eeq{7.6}
By Cauchy--Schwarz inequality, integration by parts and elementary embedding theorems we get, for arbitrary $\bfphi\in\cald_0^{1,2}(\Omega)$,
$$
|\big(\bar{(\bfsigma-\bsfw)\cdot\nabla\bsfw},\bfphi)|\le (\|\bfsigma\|_{L^4_\sharp}+\|\bfw\|_{L^4(L^4)})\|\bfw\|_{L^4(L^4)}\|\nabla\bfphi\|_2\le c\, \|\bfw\|_{\mbox{\scriptsize $\W$}}^2\|\nabla\bfphi\|_2\,.
$$

By a similar argument,
$$
\|\bar{(\bfsigma-\bsfw)\cdot\nabla\bsfw}\|_2\le (\|\bfsigma\|_{L^4_\sharp}+\|\bfw\|_{L^4(L^4)})\|\nabla\bfw\|_{L^4(L^4)} \le c\, \|\bfw\|_{\mbox{\scriptsize $\W$}}^2\,,
$$
which completes the proof of \eqref{7.6}.

We next observe that, by \theoref{exi} we have $\bsfu_0,\tilde{\bfu_0}\in W^{1,\infty}(\Omega)$, whereas $\bsfu\in L^\infty(\Omega)$, by the continuous embedding of $X^2$ into $L^\infty$. Thus, from \eqref{7.4} and also bearing in mind  the proof of the first property in \eqref{7.5} just given, one realizes that in order to show the second property in \eqref{7.5}, it is enough to prove $(\bfsigma-\bsfw)\cdot\nabla\bsfu\in L^2(L^2)$. This, in turn, follows at once  because  $\bfsigma\in W_\sharp^1\subset L^\infty(0,2\pi)$, $\bsfw\in L^2(W^{2,2})\subset L^2(L^\infty)$ and $\bsfu\in X^2(\Omega)\subset D^{1,2}(\Omega)$.
\end{proof}\par
Set $\mathscr N_i:=\mathscr P\bfN_i$, $i=1,2$, and define $\mathscr L_1(\equiv \mathscr L_{\lambda_o})$ by
\be
\mathscr L_{\lambda_o}:\bsfU\in \mathcal X\mapsto\hat{\mathscr A}(\bsfU) +\lambda_o[\hat{\partial_1}(\bsfU)+\mathscr P\mathscr C(\bsfU)]\in \caly\,.
\eeq{2.2lin}
Recalling the definition of the operator $\mathscr L_2$ given in \eqref{L2-2}, we deduce that equations \eqref{7.1} and \eqref{7.3}
can equivalently be rewritten as follows
\be\ba{ll}\medskip
\mathscr L_1(\bsfU)=\mathscr N_1(\bsfU,\bsfW,\mu)\ \ \mbox{in $\mathcal Y$}\\
\zeta\partial_\tau\bsfW+\mathscr L_2(\bsfW)=\mathscr N_2(\bsfU,\bsfW,\mu)\ \ \mbox{in ${\sf L}^2_\sharp$} \,,
\ea
\eeq{000}
where $\bsfU:=(\bsfu,\bar{\bfeta})\in X^2(\Omega)\times\real^3$ and $\bsfW:=(\bsfw,\bfxi)\in \W\times W^1_\sharp$.

\subsection{Appropriate Formulation of the Assumptions}\label{sec:Bif}
Our final goal is to describe necessary and sufficient conditions for time-periodic bifurcation. The latter means that problem \eqref{000} possesses a non-trivial family of solutions $(\bsfU(\mu),\bsfW(\mu),\zeta(\mu))$, $\mu$ in a neighborhood of 0, such that
\be\ba{ll}\medskip
{\rm (i)}\ \ \bsfW(\mu) \ \ \mbox{is $2\pi$-periodic}\,;\\
{\rm (ii)}\ \
\text{For some }\zeta_0\neq 0\,,\, \ (\bsfU(\mu),\bsfW(\mu),\zeta(\mu))\to (\0,\0,\zeta_0), \text{ as }\mu\to 0,\\ \text{ in the corresponding spaces}\,.
\ea
\eeq{Sch1}
To this end, we begin to consider, separately, the assumptions (H1)--(H4) made in  \theoref{3.1_ar} and  to formulate them  appropriately for the case at hand, namely
\be\ba{ll}\medskip
\mathscr L_1(\bsfU)=\mathscr N_1(\bsfU,\bsfW,\mu)\ \ \mbox{in $\mathcal Y$}\\
\zeta\partial_\tau\bsfW+\mathscr L_2(\bsfW)=\mathscr N_2(\bsfU,\bsfW,\mu)\ \ \mbox{in ${\sf L}^2_\sharp$} \,,
\ea
\eeq{7.7}
with $\bsfU:=(\bsfu,\bar{\bfeta})\in X^2(\Omega)\times\real^3$ and $\bsfW:=(\bsfw,\bfxi)\in \W\times W^1_\sharp$.

Assumption ({H1}) requires $\mathscr L_1$  to be a homeomorphism.
It is shown in \cite[Lemma 15]{BoGaGa1} that the operator $\mathscr L_1$ is Fredholm of index 0. Therefore, the assumption (H1)  is satisfied if ${\sf N}[\mathscr L_1]=\{\0\}$, that is,
\be
\mathscr L_1(\bsfU)=\0\ \ \Longrightarrow\ \ \bsfU=\0\,.
\tag{H1$^\prime$}
\eeq{H1'}
By definition of $\mathscr L_1$, \eqref{H1'} is equivalent to the following request: If $(\bsfu,{\sf r},\bar{\bfeta})\in X^2(\Omega)\times D^{1,2}(\Omega)\times\real^3$ is a solution to the problem
$$\ba{cc}\medskip\left.\ba{r}\medskip
-\lambda_o(\partial_1\bsfu-\bsfu_0\cdot\nabla\bsfu-\bsfu\cdot\nabla\bsfu_0)-\Delta\bsfu+\nabla {\sf r}=\0\\
\Div\bsfu=0\ea\right\}\ \ \mbox{in $\Omega$\,,}\\ \medskip
\bsfu=\0\ \ \mbox{on $\partial\Omega$}\,,\\
\mathbb{A}  \cdot \bar{\bfeta}+\varpi\Int{\partial\Omega}{}\mathbb T(\bsfu,\bar{\sfp})\cdot\bfn=\0\,,
\ea
$$
then, necessarily, $\bsfu\equiv\nabla {\sf r}\equiv\bar{\bfeta}\equiv\0$. According to \cite[Theorem 17]{BoGaGa1}, this  implies that the equilibrium configuration ${\sf s}(\lambda)$ in \eqref{sfaco} is unique for all $\lambda\in U(\lambda_o)$.

According to \theoref{5.1}, the operator $\mathscr L_2-{\rm i}\,\zeta\,\cali$ is Fredholm of index 0 and, moreover $\Sigma:=\sigma(\mathscr L_2)\cap \{{\rm i}\real\}$ is constituted only by eigenvalues of finite algebraic multiplicity (a.m.). Therefore, the assertions (i)--(ii) of Assumption (H2) can be formulated as follows:
\be
\mbox{(i) there is $\nu_0:={\rm i}\zeta_0\in\Sigma$ with a.m. $1$;  (ii) $k\nu_0\not\in \Sigma$, for all $|k|>1$.}
\tag{H2$^\prime$}
\eeq{H2'}
Taking into account the definition of $\mathscr L_2$ given in \eqref{L2}, we show that \eqref{H2'} implies, in particular, that the eigenvalue problem
\be\ba{cc}\medskip\left.\ba{r}\medskip
-\i\zeta_0\bfw-\lambda_o\,[\partial_1\bfw-\bfu_0\cdot\nabla\bfw+(-i\zeta_0{\bfxi}-\bfw)\cdot\nabla\bfu_0]=\Delta\bfw-\nabla {\sfp}\\
\Div\bfw=0\ea\right\}\ \ \mbox{in $\Omega$}\,,\\ \medskip
\bfw(x)=-\i\zeta_0{{\bfxi}}\,, \ \mbox{ $x\in\partial\Omega$}\,,
\\
(-\zeta_0^2+\mathbb{A}  \cdot )\bfxi+\varpi\Int{\partial\Omega}{} \mathbb T(\bfw,{\sfp})\cdot\bfn=\0\,,
\ea
\eeq{8.1}
has a corresponding one-dimensional eigenspace $(\bfw,\bfxi)\in  Z^{2,2}\times\real^3$.

In view of \theoref{3.1}, the operator
$$
\mathscr Q:\bsfW\mapsto \zeta_0\partial_t\bsfW+\mathscr L_2(\bsfW)\,,
$$
is Fredholm of index $0$ so that ({H3}) is automatically satisfied in our case.

As observed in the discussion of Assumption ({H1}), its formulation \eqref{H1'} implies the uniqueness of the velocity field $\bfu_0(\mu+\lambda_o)$  of the solution  determined in \theoref{exi} and corresponding to $\mu+\lambda_o$.
Even more, the assumption \eqref{H1'} combined with \cite[Theorem 17]{BoGaGa1} entails that the map
$$
\mu\in U(0)\subset\real\mapsto \bfu_0(\mu+\lambda_o)\in X^2(\Omega)
$$
is analytic. In addition,
the nonlinear operators $\mathscr N_i$, $i = 1, 2,$ are
(at most) quadratic in $(\bsfu,\bsfw)$  and then, by \lemmref{7.1}, analytic in those variables. Therefore, we conclude that also assumption (H4) is satisfied in the case at hand.\par
Before stating the bifurcation theorems in the next subsection, our final comment regards the assumption \eqref{nupr} and its formulation in the context of our problem. To this end, we deduce from \eqref{7.4} that, in the case at hand, we have
$$
S_{011}=\partial_1\bsfw-\bsfu_0\cdot\nabla\bsfw+(\bfsigma-\bsfw)\cdot\nabla\bsfu_0
+\lambda_o\Big[\tilde{\bfu}^\prime_0(0)\cdot\nabla\bsfw+(\bsfw-\bfsigma)\cdot\nabla\tilde{\bfu}_0^\prime(0)\Big]
$$
where the prime denotes differentiation with respect to $\mu$. Therefore, denoting by $\nu = \nu(\mu)$ the
eigenvalues of $\mathscr L_2 + \mu S_{011}$, we may apply \cite[Proposition 79.15 and Corollary 79.16]{Z1} to show  that  the map $$\mu\in U(0)\mapsto \nu(\mu)\in\comp$$ is well defined and of class $C^\infty$.

\subsection{Necessary and Sufficient Conditions for a Time-Periodic Bifurcation}

We are now in a position to state our bifurcation results. We begin with a necessary condition.
\Bt
Suppose there exists $(\lambda_o,\bsfu_0)$ such that   \eqref{H1'} hold. Necessary condition for the occurrence of time-periodic bifurcation is that ${\rm i}\zeta_0\in \sigma(\mathscr L_2)$.
\ET{NEC}
\begin{proof}
Letting $\beta:=\zeta_0-\zeta$, \eqref{7.7} furnishes
\be
\mathscr L_1(\bsfU)=\mathscr N_1(\bsfU,\bsfW,\mu)\mbox{ in $\mathcal Y$,}\qquad
\zeta_0\partial_\tau\bsfW+\mathscr L_2(\bsfW)=\beta\,\partial_\tau\bsfW+\mathscr N_2(\bsfU,\bsfW,\mu)\ \ \mbox{in ${\sf L}^2_\sharp$}\,.
\eeq{7.711}
Clearly, $\bsfU=\bsfW=\0$ is a solution to \eqref{7.711}, for all $(\beta,\mu)\in \real\times\real$.
By \eqref{H1'} $\mathscr L_1$ is a homeomorphism.
Assume that ${\rm i}\zeta_0\not\in \sigma(\mathscr L_2)$.
Then, by \theoref{3.1} and \theoref{5.1}, $\zeta_0\partial_\tau+\mathscr L_2$ is a homeomorphism  and,  therefore, by the implicit function theorem, for $(\beta,\mu)$ in a neighborhood of $(0,0)$, $\bsfU=\bsfW=\0$ is the only solution, and bifurcation does not occur.
\end{proof}

The following result provides sufficient conditions and is a direct consequence of  \theoref{3.1_ar}.
\Bt
Suppose there exists $(\lambda_o,\bsfu_0,\zeta_0)$ such that  assumptions \eqref{H1'} and \eqref{H2'} hold and, moreover,
$$
\Re[\nu^\prime(0)]\neq 0\,.
$$
Then, the following properties are valid. \smallskip\\
{\rm (a)} {\rm Existence.} There is an analytic family
\be
\big({\bsfu}(\varepsilon),\bsfw(\varepsilon),\bfxi(\varepsilon),\zeta(\varepsilon),\mu(\varepsilon)\big)\in X^2(\Omega)\times {\sf W}^2_\sharp \times W^1_\sharp\times \real_+\times\real
\eeq{fam1}
of solutions to \eqref{7.1}, \eqref{7.3}, $\varepsilon$ in a neighborhood $\mathcal I(0)$ of\, $0\in\real$. Moreover, let $(\bsfw_0,\bfxi_0)\in Z^{2,2}_{\mbox{\tiny $\comp$}}\times \comp^3$ be a normalized eigenfunction of the operator $\mathscr L_2$ corresponding to the eigenvalue $\i\zeta_0$, and set $(\bsfw_1,\bfxi_1):=\Re[(\bsfw_0,\bfxi_0)\,{\rm e}^{-\i\tau}]$. Then
\be
\big({\bsfu}(\varepsilon),\bsfw(\varepsilon)-\varepsilon\,\bsfw_1,\bfxi(\varepsilon)-\varepsilon\,\bfxi_1,\zeta(\varepsilon),\mu(\varepsilon)\big)\to (0,0,\zeta_0,0)\ \ \mbox{as $\varepsilon\to 0$}\,.
\eeq{Ar.101}
\par\noindent
{\rm (b)} {\rm Uniqueness.}
There is a neighborhood  $$\calu(\0,\0,\0,\zeta_0,0)\subset X^2(\Omega)\times {\sf W}^2_\sharp\times W^1_\sharp\times \real_+\times \real$$ such that every (nontrivial) $2\pi$-periodic solution to \eqref{7.1}, \eqref{7.2},  lying in $\calu$ must coincide, up to a phase shift, with a member of the family \eqref{fam1}.
\smallskip\par\noindent
{\rm (c)} {\rm Parity.}  The functions $\zeta(\varepsilon)$ and $\mu(\varepsilon)$ are even:
$$
\zeta(\varepsilon)=\zeta(-\varepsilon)\,,\ \ \mu(\varepsilon)=\mu(-\varepsilon)\,,\ \ \mbox{for all $\varepsilon\in\cali(0)$\,.}
$$
\ET{8.1}
The bifurcation due to these solutions is either subcritical or supercritical, a two-sided bifurcation being excluded by  (c), unless $\mu\equiv 0$.
\appendix

\renewcommand{\theequation}{A.\arabic{equation}}

\section{Proof of \theoref{3.1_ar}}\label{sec:proof4.1}
We proceed in several steps,  beginning with some preparatory results.
Let $v_0$ be a normalized eigenvector of $L_2$ corresponding to the eigenvalue $\nu_0$, and set
\be
v_1:=\Re[v_0\,{\rm e}^{-{\rm i}\,\tau}]\,,\ \ v_2:=\Im[v_0\,{\rm e}^{-{\rm i}\,\tau}]\,.
\eeq{vee}
Then we prove:
\Bl {\sl Let $Q$ be as in} {\rm (H3)}. {\sl Under the assumption {\rm (H2)}, we have ${\rm dim}\,{\sf N}\,[Q]=2$, and $\{v_1,v_2\}$ is a basis in ${\sf N}\,[Q]$.}
\label{3.1_00}
\EL{3.1_00}
\begin{proof} Clearly, $\mathcal S:= {\rm span}\,\{v_1,v_2\}\subseteq {\sf N}[Q]$. Conversely, take $w\in {\sf N}[Q]$, and expand it in Fourier series
$$
w=\sum_{\ell=-\infty}^\infty w_\ell\,{\rm e}^{-{\rm i}\,\ell\,\tau}\,;\ w_\ell:=\frac1{2\pi}\int_{0}^{2\pi} w(\tau)\,{\rm e}^{{\rm i}\,\ell\,\tau}\,d\tau\,,\ \ w_0\equiv\bar w=0.
$$
Obviously, $w_\ell\in \calw_{\comps}\equiv {\sf D}_{\comps}[L_2]$. From $Q(w)=0$ we deduce
$$
-\ell\,\nu_0\,w_\ell+L_2(w_\ell)=0\,,\ \ w_\ell\in{\sf D}_{\comps}[L_2]\,,\ \ \ell\in\mathbb Z,
$$
which, by (H2) and the fact that $w_0=0$, implies $w_\ell=0$ for all $\ell\in \mathbb Z\backslash\{\pm 1\}$. Thus, recalling that $\nu_0$ is simple, we infer  $w\in\cals$ and the lemma follows.\end{proof}
\smallskip\par
Denote by $\langle\cdot,\cdot\rangle$ the scalar product in $\calz$ and set
$$
(w_1|w_2):=\int_{0}^{2\pi}\langle w_1(s),w_2(s)\rangle\,{\rm d}s\,,\ \ w_1,w_2\in \calz_{2\pi,0}\,.
$$
Let $L^\dagger_2$ be the adjoint of  $L_2$. Since $\nu_0$ is simple and $L_2-\nu_0I$ is Fredholm of index 0 (by (H2)), from classical results (see e.g.\ \cite[Section 8.4]{Z}), it follows that there exists at least one  element  $v_0^\dagger\in{\sf N}_{\comps}[L_2^\dagger-\nu_0\,I]$ such that $\langle v_0^\dagger,v_0\rangle\neq 0$. Without loss, we may take
\be
\langle v_0^\dagger,v_0\rangle=\pi^{-1}\,.
\eeq{3.6_00}

We then define
$$
v_1^\dagger:=\Re[v_0^\dagger\,{\rm e}^{{\rm i}\,\tau}]\,,\ \ v_2^\dagger:=\Im[v_0^\dagger\,{\rm e}^{{\rm i}\,\tau}]\,,
$$
and observe that,
by \eqref{vee} and \eqref{3.6_00},
\be
(v_1|v_1^\dagger)=(v_2|v_2^\dagger)=1\,,\ \ (v_2|v_1^\dagger)=(v_1|v_2^\dagger)=0\,,\qquad
( (v_1)_\tau|v_1^\dagger)=0\,,\ \ ((v_1)_\tau|v_2^\dagger)=-1\,.
\eeq{3.7_00}
Set
$$
\hat{\calz}_{2\pi,0}=\big\{w\in {\calz}_{2\pi,0}: \ (w|v_1^\dagger)=(w|v_2^\dagger)=0\big\}\,,\ \
\hat{\calw}_{2\pi,0}={\calw}_{2\pi,0}\cap \hat{\calz}_{2\pi,0}\,.
$$

Let us now show that $Q$ is a homeomorphism.

\Bl {\sl Let {\rm (H2)} and {\rm (H3)} hold. Then, the operator $Q$ maps  $\hat{\calw}_{2\pi,0}$ onto $\hat{\calh}_{2\pi,0}$ homeomorphically.}
\label{3.2_00}
\EL{3.2_00}
\begin{proof} By (H3), $Q$ is Fredholm of index 0, whereas, by Lemma \ref{3.1_00}, we know ${\rm dim}\,{\sf N}\,[Q]=2$. From  classical theory of Fredholm operators (e.g. \cite[Proposition 8.14(4)]{Z}) it then follows that ${\rm dim}\,{\sf N}\,[Q^\dagger]=2$
where
$$
Q^\dagger=\zeta_0(\cdot)_\tau+L_2^\dagger
$$
is the adjoint of $Q$. In view of the stated properties of $v_0^\dagger$, we infer ${\rm span}\,\{v_1^\dagger,v_2^\dagger\}={\sf N}\,[Q^\dagger]$, and the lemma follows from another classical result on Fredholm operators (e.g. \cite[Proposition 8.14(2)]{Z}).\end{proof}
\smallskip\par
Let
$$
L_2(\mu):=L_2+\mu\,S_{011}\,.
$$
since, by (H2), $\nu_0$ is a simple eigenvalue of $L_2(0)\equiv L_2$,  denoting by $\nu(\mu)$ the eigenvalues of $L_2(\mu)$, it follows (e.g.  \cite[Proposition 79.15 and Corollary 79.16]{Z1}) that in a neighborhood of $\mu=0$ the map $\mu\mapsto\nu(\mu)$ is well defined and of class $C^\infty$ and, further,
$$
\nu'(0)=\langle v_0^\dagger, S_{011}(v_0)\rangle
.
$$
Using the latter, by  direct inspection we show
\be
\Re[\nu'(0)]=\pi^{-1}(S_{011}(v_1)|v_1^\dagger)\,.
\eeq{nu0}
\par\medskip
\begin{proof}[Proof of Assertion \rm (a)]
In order to ensure that the solutions we are looking for are non-trivial, we endow \eqref{Ar.5} with the side condition
\be
(w|v_1^\dagger)=\varepsilon\,,\ \ (w|v_1^\dagger)=0\,,
\eeq{Ar.8}
where $\varepsilon$ is a real parameter ranging in a neighborhood, $\cali(0)$, of $0$.
We next scale $v$ and $w$  by setting
$v=\varepsilon\,{\sf v}$, ${w}=\varepsilon\, {\sf w}$, so that problem \eqref{Ar.5},
 \eqref{Ar.8} becomes
\be\ba{ll}\medskip
L_1({\sf v})={\mathcal N}_1(\varepsilon, {\sf v},{\sf w},\mu)\,,\ \mbox{in $\calv$}\,;\\
\zeta_0\, {\sf w}_\tau +L_2({\sf w)={\mathcal N}_2(\varepsilon, \zeta, {\sf v},{\sf w},\mu)\,,\ \mbox{in $\calz_{2\pi,0}$}}\,,
\ \
({\sf w}|v_1^\dagger)=1\,,\ \ ({\sf w}|v_1^\dagger)=0\,,
\ea
\eeq{3.11}
where
$$\ba{ll}\medskip
\mathcal N_1(\varepsilon,\sfv,\sfw,\mu):=(1/\varepsilon)\,N_1(\varepsilon{\sf v}, \varepsilon\sfw,\mu)\,, \\
\mathcal N_2(\varepsilon,\zeta,\sfv,\sfw,\mu):=(1/\varepsilon)\,N_2(\varepsilon{\sf v}, \varepsilon\sfw,\mu)+(\zeta_0-\omega){\sf w}_\tau\,.\ea
$$ Set ${\sf U}:=(\mu, \zeta,{\sf v},{\sf w})$, and consider the map
$$\ba{cc}\medskip
F:\cali(0)\times \left(U(0)\times V(\zeta_0)\times \calu\times \calw_{2\pi,0}\right)\to \calv\times \calz_{2\pi,0}\times \real^2\,,\\ \smallskip
(\varepsilon,{\sf U})
\mapsto
\Big( L_1(\sfv)-\mathcal N_1(\varepsilon,\sfv,\sfw,\mu),\
Q(\sfw)-\mathcal N_2(\varepsilon,\zeta,\sfv,\sfw,\mu),\
({\sf w}|v_1^\dagger)-1,\ ({\sf w}|v_2^\dagger)\Big),
\ea
$$
with $U(0)$ and $V(\zeta_0)$  neighborhoods of 0 and  $\zeta_0$, respectively. Since, by (H4), we have in particular $\mathcal N_1(0,0,v_1,0)=
\mathcal N_2(0,\zeta_0,v_1,0)=0$, using
\eqref{3.7_00}$_1$ and Lemma \ref{3.1_00} we deduce that, at $\varepsilon=0$,
the equation $F(\varepsilon,{\sf U})=0$ has the solution ${\sf U}_0=(0,\zeta_0,0,v_1)$. Therefore, since by (H4) we have that $F$ is analytic at $(0,{\sf U_0})$, by the  analytic version of the Implicit Function Theorem (e.g. \cite[Proposition 8.11]{Z}), to show the existence statement --including the validity of \eqref{Ar.10}-- it suffices to show that the Fr\'echet derivative, $D_{{\mbox{\tiny {\sf U}}}}F(0,{\sf U}_0)$ is a bijection.
Now,   in view of the assumption (H4), one can easily check that the Fr\'echet derivative of $\mathcal N_1$ at $(\varepsilon=0, \sfv=0,\sfw=v_1,\mu=0)$ is equal to 0, while that of $\mathcal N_2$ at $(\varepsilon=0, {\sf U}={\sf U}_0)$ is equal to $ -\zeta\,(v_{1})_\tau+\mu\,S_{011}(v_1)$\,.
Therefore, $D_{{\mbox{\tiny {\sf U}}}}F(0,{\sf U}_0)$ is a bijection  if  we  prove that for any
$({\sf f}_1,{\sf f}_2,{\sf f}_3,{\sf f}_4)\in \calv\times\calz_{2\pi,0}\times\real\times\real$, the following system of equations has one and only one solution $(\mu,\zeta,\sfv,\sfw)\in \real\times\real\times \calu\times \calw_{2\pi,0}$:
\be\ba{rl}\medskip
{L}_1(\sfv)=&\!\!\!\!{\sf f}_1\ \ \mbox{in $\calv$}\\ \medskip
Q(\sfw)= &\!\!\!\!-\zeta\, (v_{1})_\tau+\mu\,S_{011}(v_1)+{\sf f}_2\ \ \mbox{in $\calz_{2\pi,0}$}\,,\\
({\sf w}|\bfv_1^\dagger)=&\!\!\!\!{\sf f}_3\,,\ \ ({\sf w}|\bfv_2^\dagger)={\sf f}_4 \ \ \mbox{in $\real$}\,,
\ea
\eeq{3.12}
In view of (H1), for any given ${\sf f}_1\in \calv$, equation \eqref{3.12}$_1$ has one and only one solution $\sfv\in \calu$. Therefore, it remains to prove the existence and uniqueness property only for the system of equations \eqref{3.12}$_{2-4}$
To this aim, we observe that, by Lemma \ref{3.2_00}, for a given ${\sf f}_2\in \calz_{2\pi,0}$,  equation \eqref{3.12}$_2$ possesses a unique solution $\sfw_1\in\hat{\calw}_{2\pi,0}$ if and only if its right-hand side is in $\hat{\calz}_{2\pi,0}$, namely,
$$
\big(-\zeta\, (v_{1})_\tau+\mu\,S_{011}(v_1)+{\sf f}_2|v_1^\dagger\big)=\big(-\zeta\,( v_{1})_\tau+\mu\,S_{011}(v_1)+{\sf f}_2|v_2^\dagger\big)=0\,.
$$
Taking into account \eqref{3.7_00}$_{2}$ the above conditions will be satisfied provided we can find $\mu$ and $\zeta$  satisfying the following algebraic system
\be
\mu(\,S_{011}(v_1)|v_1^\dagger)=-({\sf f}_2|v_1^\dagger)\,,\qquad
\zeta+\mu\,(\,S_{011}(v_1)|v_2^\dagger)=-({\sf f}_2|v_2^\dagger)\,.
\eeq{3.13}

However, by virtue of \eqref{nu0}, \eqref{nupr} this system possesses a uniquely determined solution $(\mu,\zeta)$, which ensures the existence of a unique solution ${\sf w}_1\in \hat{\calw}_{2\pi,0}$ to \eqref{3.12}$_2$ corresponding to the  selected values of $\mu$ and $\zeta$.\par
We now set
$$
\sfw:={\sf w}_1+\alpha\,v_1+\beta\,v_2\,,\ \ \alpha\,,\, \beta\in\real\,.
$$
Clearly, by Lemma \ref{3.1_00}, ${\sfw}$ is also a solution to \eqref{3.12}$_2$. We then
choose $\alpha$ and $\beta$ in such a way that $\sfw$ satisfies both conditions \eqref{3.12}$_{3,4}$  for any given ${\sf f}_i\in\real$, $i=1,2$. This choice is made possible by virtue of \eqref{3.7_00}$_1$.
We have thus shown that $D_{\mbox{\tiny {\sf U}}}F(0,{\sf U}_0)$ is surjective.\par
To show that it is also injective,
set ${\sf f}_i=0$ in \eqref{3.12}$_{2-4}$. From \eqref{3.13} and \eqref{nu0}, \eqref{nupr} it then follows $\mu=\zeta=0$ which in turn implies, by \eqref{3.12}$_2$ and  Lemma \ref{3.1_00}, $\sfw=\gamma_1\,v_1+\gamma_2\,v_2$, for some $\gamma_i\in\real$, $i=1,2$. Replacing this information back in \eqref{3.12}$_{3,4}$ with ${\sf f}_3={\sf f}_4=0$, and using \eqref{3.7_00}$_1$ we conclude $\gamma_1=\gamma_2=0$, which proves  the claimed injectivity property.
Thus, $D_{\mbox{\tiny {\sf U}}}F(0,{\sf U}_0)$ is a bijection, and the proof of the existence statement in (a) is completed.
 \end{proof}

\begin{proof}[Proof of Assertion \rm (b)]
Let $(z,s)\in \calu\times\calw_{2\pi,0}$  be a  $2\pi$-periodic solution to \eqref{Ar.5} with $\zeta\equiv\tilde{\zeta}$ and $\mu\equiv\tilde\mu$.
By the uniqueness property associated with  the implicit function theorem, the proof of the claimed uniqueness
amounts to show that we can find a sufficiently small $\rho>0$ such that if
\be
\|z\|_{\calu}+\|s\|_{\calw_{2\pi,0}}+|\tilde\zeta-\zeta_0|+|\tilde\mu|<\rho\,,
\eeq{3.14}
then there exists a neighborhood of $0$, $\cali(0)\subset\real$, such that
\be\ba{cc}\medskip
s=\eta\, v_1+\eta\,{\sf s}\,,\, \ z=\eta\,{\sf z}\,, \ \mbox{for all $\eta\in\cali(0)$},\, \\
|\tilde\zeta-\zeta_0|+|\tilde\mu|+\|{\sf z}\|_{\calu}+\|{\sf s}\|_{\calw_{2\pi,0}}\to 0\ \ \mbox{as $\eta\to 0$}\,.\ea
\eeq{3.15}
To this end, we notice that, by \eqref{3.7_00}$_1$, we may write
\be
s={\sigma} +\ts
\eeq{3.16}
where ${\sigma}=(s|v^\dagger_1)\,v_1+(s|v^\dagger_2)\,v_2$ and
\[
(\tilde{\sf s}|v^\dagger_i)=0\,,\ \ i=1,2\,.
\]

We next make the simple but important observation that if we modify $s$ by a  constant phase shift in time, $\delta$, namely, $s(\tau)\to s(\tau+\delta)$,  the shifted function is still  a $2\pi$-periodic solution to \eqref{Ar.5}$_2$ and, moreover, by an appropriate choice of $\delta$,
\be
{\sigma}=\eta\, v_1\,,
\eeq{3.18-A}
with $\eta=\eta(\delta)\in\real$. (The proof of \eqref{3.18-A} is straightforward, once we take into account the definition of $v_1$ and $v_2$.)
Notice that from \eqref{3.14}, \eqref{3.16}--\eqref{3.18-A} it follows that
\be
|\eta| +
\|{\ts}\|_{\calw_{2\pi,0}}\to 0 \ \ \mbox{as $\rho\to 0$}\,.
\eeq{3.19-A}
From \eqref{3.5} we thus get
\be
{L}_1(z)= N_1(z,\eta\, v_1+\ts,\tilde\mu)
\eeq{3.20}
and, recalling Lemma \ref{3.1_00},
\be
Q(\ts)=\eta(\zeta_0-\zeta)(v_1)_\tau+(\zeta_0-\zeta)\ts_\tau+N_2(z,\eta\,v_1+\ts,\tilde\mu)\,.
\eeq{3.21}

In view of (H4) and \eqref{3.14}, we easily deduce
$$
N_1(z,\eta\, v_1+\ts,\tilde\mu)=R_{110}z(\eta\,v_1+\ts)+R_{101}z\tilde\mu+R_{020}(\eta\,v_1+\ts)^2+n_1(z,\eta,\ts,\tilde\mu)\,,
$$
where
$$
\|n_1(z,\eta,\ts,\tilde\mu)\|_{\calv}\le \epsilon(\rho) \,\left(\|z\|_{\calu}+\|\ts\|_{\calw_{2\pi,0}}+\eta^2\right)\,,\ \ \
\epsilon (\rho)\to 0\ \mbox{as}\ \rho\to 0\,,$$
so that, by \eqref{3.20} and (H1), by taking $\rho$ sufficiently small we obtain
\be
\|z\|_{\calu}\le c_1\,\big(|\eta|^2+\|\ts\|_{\calw_{2\pi,0}}^2+\epsilon(\rho)\|\ts\|_{\calw_{2\pi,0}}\big)\,.
\eeq{3.22}
Likewise,
\be\ba{rl}\medskip
N_2(z,\eta\,v_1+\ts,\tilde\mu)
=&\!\!\!\!S_{011}(\eta\,v_1+\ts)\tilde\mu+S_{110}z(\eta\,v_1+\ts)+S_{101}z\tilde\mu\\&+S_{200}z^2+S_{020}(\eta\,v_1+\ts)^2+n_2(z,\eta,\ts,\tilde\mu)\,,
\ea
\eeq{3.23}
where $n_2$ enjoys the same property as $n_1$.
From \eqref{3.21},  \eqref{3.23} and \eqref{3.7_00}$_1$ we infer, according to Lemma \ref{3.2_00}, that the
following (compatibility) conditions must be satisfied
$$\ba{ll}\medskip
-\eta\,\tilde\mu\,(S_{011}(v_1)|v_1^\dagger)=\big((\zeta_0-\zeta)\ts_\tau+ S_{011}\ts\tilde\mu+S_{110}z(\eta\,v_1+\ts)|v_1^\dagger\big)\\ \medskip
\hspace*{3.5cm}+\big(S_{200}z^2+S_{020}(\eta\,v_1+\ts)^2|v_1^\dagger\big)+(n_2|v_1^\dagger)\\ \medskip
\eta\,(\zeta-\zeta_0)=\big((\zeta_0-\zeta)\ts_\tau+ S_{011}\ts\tilde\mu+S_{110}z(\eta\,v_1+\ts)|v_2^\dagger\big)\\ \medskip
\hspace*{3.5cm}+\big(S_{200}z^2+S_{020}(\eta\,v_2+\ts)^2|v_2^\dagger\big)+(n_2|v_2^\dagger)\,,
\ea
$$
so that, from \eqref{3.9} and the property of $n_2$  we show
\be \ba{ll}\medskip
|\eta|\,\big(|\tilde\mu|+|\zeta-\zeta_0|\big)\le c_2 \big(|\zeta-\zeta_0|+|\tilde\mu|\big)\,\|\ts\|_{\calw_{2\pi,0}}
+|\eta|\,\|z\|_{\calu}+\|z\|_{\calu}^2\\
\hspace*{3.5cm}+\|\ts\|_{\calw_{2\pi,0}}^2+\eta^2\big)+\epsilon(\rho)\big(\|z\|_{\calh}+\|\ts\|_{\calw_{2\pi,0}}\big)\,.\ea
\eeq{3.24}

Finally, applying Lemma \ref{3.2_00} to \eqref{3.21} and using \eqref{3.23}, \eqref{3.14} with $\rho$ sufficiently small we get
\be
\|\ts\|_{\calw_{2\pi,0}}\le c_3\,\big(|\eta|\,(|\tilde\mu|+|\zeta-\zeta_0|)+(|\eta|+|\tilde\mu|+\epsilon(\rho))\,\|z\|_{\calu}+\|z\|_{\calu}^2+\eta^2\big)\,.
\eeq{3.25}
Summing side by side \eqref{3.22}, \eqref{3.24} and $(1/(2c_3))\times$\eqref{3.25}, and taking again $\rho$ small enough, we thus arrive at
$$
|\eta|\,\big(|\tilde\mu|+|\zeta-\zeta_0|\big)+\|z\|_{\calu}+\|\ts\|_{\calw_{2\pi,0}}\le c_4\,\eta^2\,,
$$
from which we infer the validity of \eqref{3.15}$_2$, thus proving the uniqueness property (b).\end{proof}

\begin{proof}[Proof of Assertion \rm (c)]
We notice that if $\big(v(-\varepsilon),w(-\varepsilon;\tau)\big)$ is the solution corresponding to $-\varepsilon$, we have $\big(w(-\varepsilon;\tau+\pi)|v_1^\dagger\big)=\varepsilon \,v_1$, which, by part (b), implies that $\big(v(-\varepsilon),w(-\varepsilon;\tau)\big)=~\big(v(\varepsilon),w(\varepsilon;\tau)\big)$  up to a phase shift. This, in turn, furnishes $\zeta(-\varepsilon)=\zeta(\varepsilon)$ and $\mu(-\varepsilon)=\mu(\varepsilon)$. From the latter and the analyticity of $\mu$ we then obtain
that either $\mu\equiv0$ or else there is an integer $k\ge 1$ and corresponding $\mu_k\in\real\backslash\{0\}$, such that
$$
\mu(\varepsilon)=\varepsilon^{2k}\mu_k+O(\varepsilon^{2k+2})\,.
$$
Thus, $\mu(\varepsilon)<0$ or $\mu(\varepsilon)>0$, according to whether $\mu_k$ is negative or positive. The theorem is completely proved.
\end{proof}

\section*{Acknowledgments}
Denis Bonheure is supported by the ARC Advanced 2020-25 ``PDEs in interaction'' and by the Francqui Foundation as Francqui Research Professor 2021-24. Giovanni P. Galdi is  partially supported by the National Science Foundation Grant DMS--2307811. Filippo Gazzola is supported by Dipartimento di Eccellenza 2023-27 of MUR (Italy), PRIN project {\it Partial differential equations and related geometric-functional inequalities}, and by INdAM.

\section*{Declarations}
\begin{itemize}
\item The authors have no relevant financial or non-financial interests to disclose.
\item The authors have no competing interests to declare that are relevant to the content of this article.
\item Data sharing not applicable to this article as no datasets were generated or analysed during the current study.
\end{itemize}

\ed
\begin{thebibliography}{99}

\bibitem{ammann}Ammann, O.H., von K\'arm\'an, T., Woodruff, G.B., {\it The failure of the Tacoma Narrows Bridge},
Technical Report, Federal Works Agency, Washington D.C. (1941)

\bibitem{argako} Arioli, G., Gazzola, F., Koch, H., Uniqueness and bifurcation branches for planar steady Navier-Stokes equations under
Navier boundary conditions, {\it J. Math. Fluid. Mech.} {\bf 23}, No.3, Paper No. 49, 20 pp. (2021)

\bibitem{arioli} Arioli, G., Koch, H., A Hopf bifurcation in the planar Navier-Stokes equations,
{\it J. Math. Fluid Mech.} {\bf 23}, No.3, Paper No. 70, 14 pp. (2021)

\bibitem{Bab0}Babenko, K.I., On the spectrum of a linearized problem on the flow of a viscous
incompressible fluid around a body (Russian). {\it Dokl. Akad. Nauk SSSR}, {\bf 262}, 64-68 (1982)

\bibitem{Bab}Babenko, K.I., Periodic solutions of a problem of the flow of a viscous
fluid around a body, {\it Soviet Math. Dokl.} {\bf 25}, 211-216 (1982)

\bibitem{BBGGP} Berchio, E., Bonheure, D., Galdi, G.P., Gazzola, F., Perotto, S., Equilibrium configurations of a symmetric body immersed
in a stationary Navier-Stokes flow in a planar channel, to appear in SIAM J. Math. Anal.

\bibitem{black} Blackburn, H.M., Henderson, R.D., A study of two-dimensional flow past an oscillating cylinder,
{\it J. Fluid Mech.} {\bf 385}, 255-286 (1999)

\bibitem{Bev}Blevins R.D., {\it Flow induced vibrations}, Van Nostrand Reinhold Co., New York (1990)

\bibitem{BoGa0}Bonheure, D., Galdi, G.P., Global Weak Solutions to a Time-Periodic Body-Liquid Interaction Problem, submitted (2023)

\bibitem{BGG}Bonheure, D., Galdi, G.P., Gazzola, F., Equilibrium configuration of a rectangular
obstacle immersed in a channel flow. {\it C. R. Math. Acad. Sci. Paris} {\bf 358}, 887-896 (2020);
updated version in arXiv:2004.10062v2 (2021)

\bibitem{BoGaGa1}Bonheure, D., Galdi, G.P., Gazzola, F., Stability of equilibria and bifurcations for a fluid-solid interaction problem (2024), preprint.

\bibitem{BoHiPaSpe} Bonheure, D., Hillairet, M., Patriarca, C., and Sperone, G.,
 Long-time behavior of an anisotropic rigid body interacting with a Poiseuille flow in an unbounded channel, submitted (2023)

\bibitem{CR1} Crandall, M.G., Rabinowitz, P.H., The Hopf bifurcation theorem in infinite dimensions,
{\it Arch. Ration. Mech. Anal.} {\bf  67}, 53-72 (1977)

\bibitem{diana}Diana, G., Resta, F., Belloli, M., Rocchi, D., On the vortex shedding forcing on suspension bridge deck,
{\it J. Wind Engineering and Industrial Aerodynamics}, {\bf 94(5)}, 341-363 (2006)

\bibitem{Dyr}Dyrbye, C., Hansen, S.O., {\it Wind Loads on Structures}, Wiley, New York (1997)

\bibitem{FN}Farwig, R., Neustupa, J., Spectral properties in $L^q$ of an Oseen operator modelling
fluid flow past a rotating body, {\it Tohoku Math. J.} {\bf 62}, 287--309 (2010)

\bibitem{FGQ} Fatone, L., Gervasio, P., Quarteroni, A.,
Multimodels for incompressible flows, {\it
J. Math. Fluid Mech.}  {\bf 2}, 126--150 (2000)

\bibitem{RR} Gerecht, D., Rannacher, R., Wollner,W., Computational aspects of pseudospectra
in hydrodynamic stability analysis. {\it J. Math. Fluid Mech.} {\bf 14}, 661--692 (2012)

\bibitem{Gah} Galdi, G.P., On the motion of a rigid body in a viscous liquid: A mathematical analysis with
applications, {\it Handbook of Mathematical Fluid Mechanics}, Elsevier Science, 653-791 (2002)

\bibitem{Gab}Galdi, G.P., {\it An introduction to the mathematical theory of the Navier-Stokes equations.
Steady-state problems}, Second edition. Springer Monographs in Mathematics, Springer, New York (2011)

\bibitem{GaCe}Galdi, G.P., Steady-state Navier-Stokes problem past a rotating body: geometric-functional properties and related questions., {\it Topics in mathematical fluid mechanics}, 109-197, Springer Lecture Notes in Math., {\bf 2073} (2013)

\bibitem{GaBif}Galdi, G.P., A time-periodic bifurcation theorem and its applications to Navier-Stokes flow past an obstacle, in {\it Mathematical Analysis of Viscous Incompressible Flow},
edited by T.~Hishida, R.I.M.S. Kokyuroku (Kyoto University, Japan, 2015), pp. 1-27. arXiv:1507.07903

\bibitem{GaArma}Galdi, G.P., On bifurcating time-periodic flow of a Navier-Stokes liquid past a cylinder. {\it Arch. Ration. Mech. Anal.}
{\bf 222}, 285-315 (2016)

\bibitem{GaMaH} Galdi, G.P., Kyed, M., Time-periodic solutions to the Navier-Stokes equations.
{\it Handbook of mathematical analysis in mechanics of viscous fluids}, 509-578, Springer, Cham (2018)

\bibitem{GPP}Gazzola, F., Pata, V., Patriarca, C., Attractors for a fluid-structure interaction problem
in a time-dependent phase space, J. Funct. Anal. {\bf 286}, Paper No. 110199, 56 pp. (2024)

\bibitem{GazP}Gazzola, F., Patriarca, C., An explicit threshold for the appearance of lift on the deck of a bridge. {\it J. Math. Fluid Mech.} {\bf 24}, No.1, Paper No. 9, 23 pp. (2022)

\bibitem{GazzSp}Gazzola, F., Sperone, G., Steady Navier-Stokes equations in planar domains with obstacle and explicit bounds for unique solvability, Arch. Ration. Mech. Anal. {\bf 238}, 1283-1347 (2020)

\bibitem{GG}Gohberg, I., Goldberg, S., Kaashoek, M.A., {\it Classes of linear operators: I. Operator Theory},  Advances and Applications, Vol.49, Birkh\"auser Verlag, Basel (1990)

\bibitem{Hey} Heywood, J.G., The Navier-Stokes equations: on the existence, regularity and decay of solutions, {\it Indiana Univ. Math. Journal} \textbf{29}, 639-681 (1980)

\bibitem{PPDL} Paidoussis, M., Price, S., De Langre, E., {\it Fluid-Structure Interactions: Cross-Flow-Induced Instabilities},
Cambridge University Press (2010)

\bibitem{Patri} Patriarca, C., Existence and uniqueness result for a fluid-structure-interaction evolution problem in an unbounded 2D channel, {\it NoDEA Nonlinear Differential Equations Appl.} {\bf 29}, No. 4, Paper No. 39, 38 pp. (2022)

\bibitem{Saz} Sazonov, L.I., The onset of auto-oscillations in a flow, {\it Siberian Math. J.} {\bf 35}, 1202-1209 (1994)

\bibitem{scott} Scott, R., {\it In the wake of Tacoma. Suspension bridges and the quest for aerodynamic stability}, ASCE Press (2001)

\bibitem{Will} Williamson, C.H.K., Govardhan, S., Vortex-induced Vibrations, {\it Ann. Rev. Fluid Mech}, {\bf 36}, 413-55 (2004)

\bibitem{Z} Zeidler, E., {\it Nonlinear Functional Analysis and Applications}, Vol.1, Fixed-Point Theorems, Springer-Verlag, New York (1986)

\bibitem{Z1}Zeidler, E., {\it Nonlinear Functional Analysis and Applications}, Vol.4, Application to Mathematical Physics, Springer-Verlag, New York (1988)

\end{thebibliography}
